\documentclass[11pt]{article}
\usepackage{amsfonts,amsmath,amsthm,amssymb}
\usepackage{graphicx}
\usepackage{booktabs}
\usepackage{float}
\usepackage{longtable}
\usepackage{geometry}
\usepackage[table]{xcolor}
\usepackage{url}
\usepackage[draft,defaultcolor=red]{changes}

\title{Compressed Traffic Assignment with the \\Augmented Lagrangian Method\footnote{
This work was carried out at the Fulton School of Engineering, Arizona State University, Tempe, AZ. X. Zhou, Y. Liu, and Y. Li are with School of Sustainable Engineering and the Built Environment, ASU. P. Li is with Norfolk Southern Corporation. D. Bertsekas is with School of Computing and Augmented Intelligence, ASU.}}
\author{Xuesong Zhou, Peiheng Li, Yajun Liu, Yuchao Li, and Dimitri Bertsekas}
\date{September 2026}

\begin{document}
\maketitle

\begin{abstract}
We consider large-scale traffic assignment problems and develop a path-based compression framework. In particular, we partition paths into major and minor paths according to nominal path flows and a prescribed threshold, and retain the major paths explicitly. For the minor paths, we use reference flows derived from nominal-flow information and represent deviations using a truncated singular value decomposition of the minor path-link incidence matrix. The resulting compressed formulation preserves convexity, and a simple scaling of the nominal minor flows ensures feasibility under changes in demand. To solve the resulting formulation, we use an augmented Lagrangian method with separate penalty parameters for the different constraints and adaptive penalty parameter updates.

We conduct computational studies on the Chicago Sketch, Chicago Regional, and Philadelphia networks in two stages. Using nominal demand data, we first study the selection of the compression threshold and rank and the effect of path-pool expansion. The results show substantial dimensionality reduction while maintaining high solution accuracy. We then reuse the nominal representation for perturbed-demand problems. The compressed formulation remains highly accurate and reduces computational time by approximately 58\%–77\% on the original path pools and 51\%–76\% on expanded path pools. We also compare it with simpler reduced formulations that either discard the minor paths or fix their flows at the reference values. The results show that fixing the minor-path flows at their reference values already provides high accuracy, while the low-dimensional adjustments can further improve link travel-time accuracy in some cases.

\end{abstract}

\section{Introduction}\label{sec:intro}
The traffic assignment problem is a classical problem in transportation science and network optimization. In its standard user-equilibrium formulation, it can be posed as a convex optimization problem that distributes travel demand across network paths subject to flow conservation and nonnegativity constraints. This formulation goes back to the seminal work of \cite{beckmann1956studies}, and it has motivated a large literature on both modeling and computation. In particular, the Frank-Wolfe method \cite{frank1956algorithm} became a standard approach because it avoids exhaustive path enumeration by solving shortest-path subproblems iteratively. Subsequent developments, including gradient projection, quasi-Newton, and conjugate-direction methods, have improved computational performance substantially; see, for example, \cite{jayakrishnan1994fast}, \cite{bertsekas1980class, bertsekas2011centralized}, \cite{mitradjieva2013stiff}, and the textbook treatments in \cite{sheffi1985urban, ahuja1993network, patriksson1994traffic, bertsekas1998network, boyles2020transportation}.

Despite this progress, path-based formulations remain difficult to solve at large scale because the number of feasible paths can be enormous. As network size grows, the associated path-link incidence matrix becomes very large, and the resulting optimization problem becomes increasingly expensive in both storage and computation. Existing methods alleviate this difficulty in different ways, for example, by generating paths selectively, refining active path sets, or exploiting origin-based structure. {\color{black}Our focus is different: we seek to replace a large collection of path variables by a lower-dimensional representation of their contribution to the aggregate link flows. This dimensionality challenge is the starting point of the present work.}

{\color{black}Our approach is motivated by a simple structural observation: in many realistic traffic assignment problems, only a relatively small subset of paths carries substantial flow, while many other paths make smaller contributions to the aggregate link flows. This suggests that the full path representation may contain substantial redundancy from the standpoint of optimization. We exploit this structure by partitioning the paths into two groups. The first group consists of major paths, which are retained explicitly as decision variables. The second group consists of minor paths. For these paths, we use reference flows derived from nominal-flow information and represent deviations from these reference flows in a low-dimensional subspace obtained from a truncated singular value decomposition (SVD) of the corresponding minor-path submatrix of the path-link incidence matrix.}

{\color{black}The proposed framework is aimed primarily at applications where closely related large-scale traffic assignment problems must be solved repeatedly. A nominal-demand problem can be used in an offline design stage to identify the major and minor paths, select the compression parameters, and construct the low-dimensional representation. This nominal information can then be reused as the OD demands change, without recomputing the major-minor partition or the SVD basis. Thus, our framework is not intended to replace a one-time full solve from scratch, but to exploit reusable nominal information in order to reduce the computational cost of solving related problem instances.}

Our construction yields a compressed path-based formulation with substantially fewer variables than the original problem while preserving its convex structure. When the target demand differs from the nominal demand, a simple scaling of the nominal minor-path flows ensures feasibility of the compressed formulation. Thus, the proposed method yields a well-defined convex optimization problem while preserving the basic linear constraint structure of the original traffic assignment model.

To solve the compressed problem, we use an augmented Lagrangian (AL) method with separate penalty parameters for the equality and inequality constraints. The resulting algorithm combines multiplier updates with efficient treatment of the simple bound constraints on the major path variables. It is particularly well suited to warm starts and repeated solution calculations. We thus provide both a compressed formulation and a computational method tailored to it.

{\color{black}Computational studies on the Chicago Sketch, Chicago Regional, and Philadelphia networks support this two-stage use of the framework. On the nominal-demand problems, substantial compression preserves very high solution accuracy, and the compression remains effective when the path pools are expanded with additional candidate paths. When the major–minor partition and SVD basis constructed from the nominal problem are reused under perturbed OD demands, the compressed formulation remains highly accurate while substantially reducing computational time. Comparisons with simpler reduced formulations show that fixing the minor-path flows at their reference values already provides high accuracy, while allowing low-dimensional adjustments can further improve link travel-time accuracy in several cases.}

The remainder of the paper is organized as follows. In Section~\ref{sec:background}, we review related work on traffic assignment algorithms, spectral methods in transportation, and reduced-order optimization. In Section~\ref{sec:compression_framework}, we develop the compressed formulation and establish its feasibility and approximation properties. In Section~\ref{sec:alm}, we describe the AL method and its algorithmic properties. In Section~\ref{sec:experiments}, we present the nominal-design and demand-perturbation computational studies. Section~\ref{sec:conclusion} concludes and discusses broader implications.

\section{Background and Related Work}\label{sec:background}

In this section, we review three strands of literature most relevant to our work. We first discuss path-based traffic assignment algorithms {\color{black}and approaches for managing large path sets}. We then examine spectral and rank-based ideas in transportation, emphasizing the distinction between matrix reconstruction and the optimization-oriented use of truncated SVD in our framework. Finally, we discuss the connection between our approach and the broader literature on reduced-order optimization.

\subsection{Traffic Assignment Problem and Solution Algorithms}

Path-based formulations of the traffic assignment problem have been studied extensively, but their central computational difficulty has remained the same: the number of feasible paths can be extremely large, and the associated path--link incidence matrix grows rapidly with network size. This creates substantial challenges in both storage and computation for realistic transportation networks.

A classical way to address this difficulty is to avoid full path enumeration. The Frank-Wolfe algorithm and its variants generate shortest paths iteratively and have long served as standard methods for solving traffic assignment problems. Other path-based approaches, including gradient projection, conjugate-gradient-type methods, and Newton-type methods, improve computational performance by refining active path sets and exploiting the smooth structure of the objective. {\color{black}These methods have substantially advanced the computational treatment of traffic assignment by avoiding or managing the potentially enormous set of feasible paths.}

Another line of work improves tractability by exploiting network structure. Origin-based and bush-based methods organize computations around origins or acyclic subnetworks \cite{gallager1977minimum,bertsekas1984second,bar2002origin,nie2010class}. Other approaches use greedy path management, decomposition, or operator-splitting techniques such as ADMM to improve scalability \cite{xie2018greedy, liu2023alternating, liu2024novel, liu2025admm}. Mesoscopic and cross-resolution modeling frameworks \cite{zhou2014dtalite, zhou2022meso}, together with open-source tools such as \texttt{Path4GMNS} \cite{Path4GMNS} and \texttt{AequilibraE} \cite{camargo2015aequilibrae}, have further improved the practical study of large-scale traffic assignment problems.

{\color{black}Our focus is complementary to these approaches. Rather than managing a large path set through path generation, active-set refinement, or decomposition, we construct a lower-dimensional representation of the contribution of a given collection of paths to the aggregate link flows. The resulting formulation therefore reduces the number of decision variables associated with the path representation itself.}

\subsection{Spectral Methods and Rank Structure in Transportation}

Spectral and low-rank ideas have appeared in transportation in several different contexts. Early work used SVD and principal component analysis to support dynamic traffic assignment and online calibration \cite{peeta1999generalized, prakash2017reducing}. Other studies examined the rank of the path--link incidence matrix from the standpoint of observability, sensor placement, and network reconstruction \cite{hu2009identification, zhou2010information}. In these works, rank is interpreted in the algebraic sense, namely, as the number of linearly independent rows or columns needed to represent the matrix exactly.

{\color{black}Our use of spectral structure is different. We do not use SVD to identify an exact basis for the path space, nor do we require the path-link incidence matrix to be well approximated by a low-rank matrix in the usual reconstruction sense. Instead, we use reference flows derived from nominal-flow information and use truncated SVD to construct a low-dimensional representation of deviations from these flows. The relevant question is therefore not how accurately the minor path-link incidence matrix can be reconstructed, but whether changes in the aggregate link-flow contribution of the minor paths can be represented accurately using a small number of singular directions.}

{\color{black}The singular values of the minor path-link incidence matrix may decay gradually, so a relatively large rank can be required for accurate reconstruction of the matrix itself. Our formulation does not seek such a reconstruction. Rather, reference minor-path flows are used, and the truncated SVD is used only to represent adjustments away from those flows. Consequently, the appropriate rank should be judged by its ability to represent changes relevant to the optimization problem, rather than by spectral reconstruction accuracy alone.}

Most recently, \cite{zhang2025lowrank} constructed reduced hypothetical path sets with real-valued (non-binary) path--link incidence matrices by fitting a surrogate Logit-based stochastic user equilibrium model. Their approach provides evidence of low-dimensional structure in path-based traffic models, but it differs from ours in both objective and construction: they replace the original assignment model with a surrogate one, whereas we use the original binary incidence matrix to construct a reduced representation of minor-path flow deviations while retaining the original Beckmann objective.

\subsection{Reduced-Order Optimization}

Low-rank reduction is also a common idea in large-scale optimization. In many problems, a high-dimensional model admits an effective lower-dimensional representation that preserves the components most relevant to the optimization objective while reducing computational cost. Such ideas appear in areas such as portfolio optimization \cite{markowitz1952portfolio, ross1976arbitrage, chopra1993effect}, machine learning \cite{vapnik1995nature, scholkopf1998nonlinear, rahimi2008random}, matrix approximation \cite{eckart1936approximation, golub1965svd, sirovich1987low, sarwar2000application, koren2009matrix, candes2009exact}, and aggregation in approximate dynamic programming \cite{bertsekas2012dynamic,bertsekas2019reinforcement}.

{\color{black}Our work brings this reduced-order viewpoint to path-based traffic assignment. The key idea is to use the path-link incidence matrix to construct a low-dimensional representation of deviations from reference minor-path flows. By restricting these deviations to a low-dimensional spectral subspace, we obtain a reduced-order convex program that remains closely connected to the original traffic assignment model while involving substantially fewer decision variables. }

This perspective places our work at the intersection of path-based traffic assignment and reduced-order optimization. The resulting formulation preserves convexity, while feasibility can be ensured by retaining at least one major path for each OD pair and adjusting the nominal minor-path flows when necessary. It can then be solved using the AL method described in Section~\ref{sec:alm}.

\section{Problem Formulation and Path-Based Compression}\label{sec:compression_framework}

In this section, we describe our path-based compression method. We start by introducing the classical traffic assignment problem. We then introduce our compression framework, and provide a condition under which the compressed problem remains feasible.


\subsection{Traffic Assignment Problem}
We consider a transportation network with $m$ directed links, $n$ paths, and $\ell$ origin-destination (OD) pairs. Let $A$ be the $\ell$ by $n$ OD-path incidence matrix (rows correspond to OD pairs, and  columns correspond to paths, with entries being either 1 or 0, depending on whether the column path corresponds to  the row OD pair, or not). Let $B$ be the $n$ by $m$ path-link incidence matrix (rows correspond to paths, and  columns correspond to links, with entries being either 1 or 0, depending on whether  the row path contains the column link, or not). There is a given nonnegative OD pair demand vector $d\in \Re^\ell$. There is also a given nonnegative offset link flow vector $v_0$, which can represent the known contributions to link flows from OD pairs having a unique path.

The traffic assignment problem is defined as
\begin{equation}
    \label{eq:original_ta}
    \begin{aligned}
        \min_{x,v}\quad & f(v)\\
    \text{s.t.}\quad &v=B'x+v_0,\ \ Ax=d, \ \ x\geq 0,
    \end{aligned}
\end{equation}
where $f$ is a given convex function of the link flow vector $v$ [specifically, the Beckmann objective associated with the Bureau of Public Roads (BPR) link travel-time functions in our numerical experiments], a prime denotes transposition, and the inequality $\geq$ is meant to be component-wise. 

Implicit in our formulation is that each OD pair has at least one path and that each path belongs to exactly one OD pair, so the problem is feasible, the matrix $A$ has rank $\ell$, and each column of $A$ is a column of the $\ell\times \ell$ identity matrix.    
Since the cost function $f$ is convex and the constraints are linear, problem \eqref{eq:original_ta} is a convex programming problem.

Despite being convex and feasible, problem \eqref{eq:original_ta} can be challenging to solve. This is because the values of $n$ and $\ell$ are often quite large for practical problems. We will next introduce a compressed problem, which involves fewer variables, thus facilitating its computational solution.

\subsection{Path Compression and Compressed Problem}
\label{sec:compression}

In many practical contexts, nominal path-flow information is available from historical data, prior computations, or coarse preliminary solutions, and can be used to identify the paths that are likely to be dominant. These are the major paths that have significant impact on the link flows $v$. In contrast, the remaining paths can be regarded as minor paths, whose contributions to the link flows are secondary. Our approach is to use nominal-flow information to construct reference flows for the minor paths and to represent deviations from these flows in a low-dimensional subspace obtained via SVD.

In particular, let $\tau$ be a threshold on the nominal flow. For OD pairs that have multiple corresponding paths,  we select major paths according to the following two rules:
\begin{itemize}
    \item[1)] For each OD pair, a path that has the largest nominal flow is a major path; 
    \item[2)] The remaining paths are split into major paths (those with nominal flow no smaller than $\tau$) and minor paths (those with nominal flow smaller than $\tau$).  
\end{itemize}
The first rule ensures that each OD pair retains at least one major path, a fact that will be used in Section~\ref{sec:convex_feasible} to establish feasibility of the compressed problem. The second rule controls the extent of compression through the threshold $\tau$.

Without loss of generality, let us assume that the first $s$ elements of the vector $x$ correspond to the major paths. As a result, we can partition the vector $x$ as
$$x=\begin{pmatrix}
    y\\
    w
\end{pmatrix},$$
where $y\in \Re^s$ and $w\in \Re^{n-s}$ represent the flows of the major and minor paths, respectively. Accordingly, we can partition the matrix $B$ as
$$B=\begin{pmatrix}
    B_1\\
    B_2
\end{pmatrix},$$
where the dimensions of $B_1$ and $B_2$ are $s$ by $m$ and $(n-s)$ by $m$, respectively. Let
$$B_2\approx U_r\Sigma_r V_r',$$
denote the rank-$r$ truncated SVD of $B_2$, where $U_r$ is an $(n-s)$ by $r$ matrix, and the dimensions of $\Sigma_r$ and $V_r$ are determined accordingly. The rank $r$ is a design parameter that determines the dimension of the reduced subspace used to represent variations in the minor-path flows.

{\color{black}Let $w_0$ denote the reference minor-path flow. For the nominal-demand problem, $w_0$ is the nominal minor-path flow. When the OD demand changes, the nominal minor-path flows associated with each OD pair are scaled by the same factor as the demand of that OD pair to obtain $w_0$.

By the Eckart–Young theorem, the matrix $U_r\Sigma_rV_r'$ obtained via truncated SVD is the best rank-$r$ approximation of $B_2$ in the spectral norm; see \cite[Theorem~2.4.8]{golub2013matrix}. Hence,
$$ \left\| B_2'(w-w_0) -(U_r\Sigma_rV_r')'(w-w_0) \right\| \leq \sigma_{r+1}\|w-w_0\|, $$
where $\sigma_{r+1}$ is the largest discarded singular value and the vector norm is Euclidean. Since $U_r$ consists of the first $r$ left singular vectors of $B_2$, we have the identity $ B_2' U_r = V_r\Sigma_r. $ It follows that
$$ (U_r\Sigma_rV_r')'(w-w_0) = B_2'U_rU_r'(w-w_0). $$

For any minor-path flow $w$, the vector
$$ U_rU_r'(w-w_0) $$
is the projection of $w-w_0$ onto the subspace spanned by the first $r$ left singular vectors of $B_2$. Writing
$$ z=U_r'(w-w_0), $$
the corresponding truncated-SVD approximation of the change in link flows is
$$ B_2'U_rz. $$
This motivates restricting the deviation of the minor-path flows from $w_0$ to this subspace and writing
$$ w=U_rz+w_0, $$
where $z\in\Re^r$ is treated as a decision variable; i.e.,
$$x=\begin{pmatrix}
    y\\
    U_rz+w_0
\end{pmatrix}.$$

Under this representation, the link-flow vector is \begin{align*}
    v=&B'x+v_0=(B_1'\;\;B_2')\begin{pmatrix}
    y\\
    U_rz+w_0
\end{pmatrix}+v_0\\
=&B_1'y+Dz+B_2'w_0+v_0,
\end{align*}
where $D=B_2'U_r$ is an $m$ by $r$ matrix.} Similarly, we partition the $A$ matrix as
$$A=\begin{pmatrix}
    A_1\ \    A_2
\end{pmatrix},$$
where $A_1$ is an $\ell$ by $s$ matrix, and $A_2$ is an $\ell$ by $(n-s)$ matrix. The term $Ax$ is written as
{\color{black}\begin{align*}
    Ax=&\begin{pmatrix}
    A_1\ \    A_2
\end{pmatrix}\begin{pmatrix}
    y\\
    U_rz+w_0
\end{pmatrix}\\
=&A_1y+Mz+A_2w_0,
\end{align*}}
where $M=A_2U_r$ is an $\ell$ by $r$ matrix. In addition, the constraint $x\geq 0$ is expressed as
$$y\geq 0,\quad U_rz+w_0\geq0.$$

With the new variables $y$ and $z$ introduced above, we obtain the following compressed problem:
\begin{equation}
    \label{eq:compressed_ta}
    \begin{aligned}
        \min_{y,z}\quad & \hat f(y,z)\\
    \text{s.t.}\quad & {\color{black}y\geq 0,\ \ A_1y+Mz=d-A_2w_0, \ \ U_rz+w_0\geq0,}
    \end{aligned}
\end{equation}
where {\color{black}$\hat f(y,z)=f(B_1'y+Dz+B_2'w_0+v_0)$}. We summarize the main notation in Table~\ref{tab:notation} given in Appendix~\ref{app:notation}.

\subsection{Convexity, Feasibility, and Optimality of the Compressed Problem}\label{sec:convex_feasible}

As discussed earlier, the original problem is convex and feasible. The compressed problem also preserves convexity. Indeed, since $f$ is convex and its argument in $\hat f(y,z)$ is affine in $(y,z)$, the function $\hat f$ is convex. Moreover, all the constraints in \eqref{eq:compressed_ta} are affine. Hence, the compressed problem is convex; see, e.g., \cite[Prop. 1.1.4]{bertsekas2009convex}

{\color{black}We next consider feasibility. Suppose first that the nominal path flows used to construct the reference flow are feasible for the nominal-demand problem. For the nominal-demand problem, let $y_0$ and $w_0$ denote the nominal major- and minor-path flows, respectively. Since these flows satisfy the OD demand constraints, $ A_1y_0+A_2w_0=d. $ Therefore,
$$ d-A_2w_0=A_1y_0\geq0, $$
and $(y,z)=(y_0,0)$ is feasible for the compressed problem.

When the OD demand changes, the reference minor-path flows associated with each OD pair are scaled by the same factor as the demand of that OD pair. Since the total nominal minor-path flow for an OD pair does not exceed its nominal demand, the scaled reference minor-path flow does not exceed the new demand. Hence,
$$ d-A_2w_0\geq0. $$

By construction, each OD pair has at least one major path. Since each column of $A_1$ is a column of the $\ell\times\ell$ identity matrix, a nonnegative vector $y$ satisfying
$$ A_1y=d-A_2w_0 $$
can be obtained by assigning the residual demand of each OD pair to one of its major paths. Thus, with $z=0$, the compressed problem is feasible.

We finally compare the optimal values of the original and compressed problems. Every feasible pair $(y,z)$ of the compressed problem defines a feasible path-flow vector of the original problem through
$$ x= \begin{pmatrix} y\\ U_rz+w_0 \end{pmatrix}. $$
Thus, the feasible set of the compressed problem corresponds to a subset of the feasible set of the original problem, and the optimal value of the compressed problem cannot be smaller than that of the original problem.

If the nominal path-flow vector used to construct $w_0$ is an optimal solution of the original nominal-demand problem, then setting $z=0$ makes this optimal path-flow vector feasible for the compressed problem. Hence, in this case, the original and compressed problems have the same optimal value. This conclusion holds for any rank $r$.

When the link travel-time functions are strictly increasing, the Beckmann objective is strictly convex in the link-flow vector. Hence, the optimal link-flow vector is unique, although the corresponding optimal path-flow vector need not be unique. This nonuniqueness is inherent in the original traffic assignment problem and is neither introduced nor resolved by the compression.
}

\begin{figure}[htbp]
\centering
\includegraphics[width=0.79\textwidth]{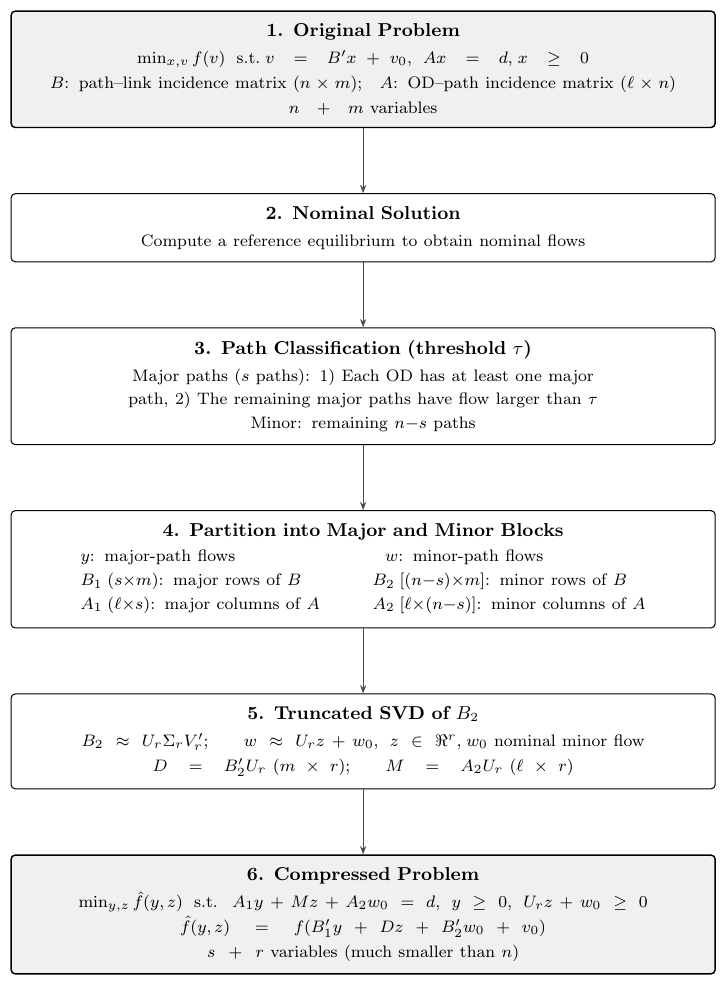}
\caption{Architecture of the major-minor path compression framework. {\color{black}A nominal solution partitions paths into major and minor paths. Major paths are retained explicitly, while deviations of the minor-path flows from their reference values are represented in a low-dimensional subspace obtained from the truncated SVD of $B_2$. The resulting formulation has $s+r$ variables, which can be substantially fewer than the $n$ path-flow variables of the original formulation.}}
\label{fig:compression-pipeline}
\end{figure}

Figure~\ref{fig:compression-pipeline} summarizes the six steps of our compression procedure, from the original problem with $n$ variables to the reduced formulation with $s + r$ variables.

\section{Solution via the Augmented Lagrangian Method}
\label{sec:alm}
{\color{black}We solve problem \eqref{eq:compressed_ta} via the AL method with partial elimination of constraints; see, e.g., \cite[Section 2.4]{Ber82a}, \cite[Section~5.2]{bertsekas2016nonlinear}. There are two reasons for this choice. First, the equality constraint and the constraint $U_rz+w_0\geq0$ can be incorporated into the augmented Lagrangian, while the simple bound constraint $y\geq0$ can be handled directly by an efficient bound-constrained optimization method. Second, the solution of one AL subproblem can be used to warm-start the next. These features make the AL method well suited for solving the compressed problem \eqref{eq:compressed_ta}. In what follows, we describe the AL function and the multiplier and penalty update formulas for each iteration.}

\subsection{Augmented Lagrangian Function}
We denote by $\lambda$ and $\mu$ the multipliers associated with the constraints {\color{black}$A_1 y + M z +A_2w_0- d = 0$ and $U_r z+w_0 \geq 0$,} respectively. As a result, we have $\lambda\in \Re^\ell$ and $\mu\in \Re^{n-s}$. For these two sets of constraints, we also associate different penalty parameters, which define a two-dimensional vector $c$, i.e., $c = (c_1, c_2)$, where $c_1$ and $c_2$ are associated with the constraints {\color{black}$A_1 y + M z+A_2w_0 - d = 0$ and $U_r z+w_0 \geq 0$}, respectively. We will also use the subscript $i$ to denote the $i$th component of a vector. In particular, {\color{black}$[U_r z+w_0]_i$ is the $i$th component of the vector $U_r z+w_0$.} 

Using the notation introduced above, the AL function is given as:
{\color{black}\begin{equation}
    \label{eq:aug_lag_form_2}
    \begin{aligned}
            {L}_c(y,z,&\lambda,\mu)\\
    =&\hat f(y,z)+\lambda'(A_1y+Mz+A_2w_0-d)+\frac{c_1}{2}\Vert A_1y+Mz+A_2w_0-d\Vert^2+\\
    &\ \ \frac{1}{2c_2}\sum_{i=1}^{n-s}\Big\{\big(\max\{0,\mu_i-c_2[U_rz+w_0]_i\}\big)^2-\mu_i^2\Big\}.
    \end{aligned}
\end{equation}}

If $f$ is continuously differentiable, as in our numerical experiments, the AL function $L_c(y, z, \lambda, \mu)$ is also continuously differentiable with respect to $y$ and $z$. The gradients are given by:
{\color{black}\begin{align*}
    \nabla_{y}{L}_c(y,z,\lambda,\mu)=&\nabla_{y}\hat f(y,z)+A_1'\big(\lambda+ c_1(A_1y+Mz+A_2w_0-d)\big),\\
    \nabla_{z}{L}_c(y,z,\lambda,\mu)=&\nabla_{z}\hat f(y,z)+M'\big(\lambda+c_1(A_1y+Mz+A_2w_0-d)\big)-
    U_r'\big(\mu+c_2h^+(z,\mu,c_2)\big),
\end{align*}}
where
{\color{black}\begin{equation}
\label{eq:h_p_def}
    h^+(z,\mu,c_2)=\begin{pmatrix}
        h_1^+(z,\mu,c_2)\\
        h_2^+(z,\mu,c_2)\\
        \vdots\\
        h_{n-s}^+(z,\mu,c_2)
    \end{pmatrix},\qquad h_i^+(z,\mu,c_2)=\max\bigg\{-[U_rz+w_0]_i,-\frac{\mu_i}{c_2}\bigg\}.
\end{equation}}
A detailed derivation of the preceding formulas for problems with general inequality constraints can be found in several textbooks, see e.g., \cite[Chapter 3]{Ber82a}.

\subsection{Iterations of the Augmented Lagrangian Method}
\label{sec:alm_iterations}
Let us discuss the updating equations for the multipliers and penalty parameters. Given the multipliers of the $k$th iteration, which we denote as $\lambda^k$ and $\mu^k$, and the penalty parameter $c^k = (c_1^k, c_2^k)$, we first compute the values $y^k$ and $z^k$ via the following minimization:
\begin{equation}
    \label{eq:inner_iteration}
    (y^k, z^k) \in \arg \min_{y \geq 0, z} L_{c^k}(y, z, \lambda^k, \mu^k).
\end{equation}
In our implementation, this inner minimization is solved using the L-BFGS-B code, which handles the simple bound constraint $y \geq 0$ by using the gradient projection method developed by the last author in \cite[Section~1.5]{Ber82a}.

Having computed the values $y^k$ and $z^k$, the updated multipliers $\lambda^{k+1}$ and $\mu^{k+1}$ are given by:
{\color{black}\begin{equation}
    \label{eq:multiplier_update}
    \begin{aligned}
        \lambda^{k+1} &= \lambda^k + c_1^k (A_1 y^k + M z^k+A_2w_0 - d), \\
    \mu^{k+1} &= \mu^k + c_2^k h^+(z^k, \mu^k, c_2^k),
    \end{aligned}
\end{equation}}
cf.\ Eq.~\eqref{eq:h_p_def}.

To accelerate convergence while controlling ill-conditioning in the AL minimization, we update separately the penalty parameters associated with different constraints. In particular, we fix some parameters $\beta > 1$ and $\gamma \in (0,1)$, and set $c^{k+1} = (c_1^{k+1}, c_2^{k+1})$ according to:
{\color{black}\begin{equation}
    \label{eq:penalty_update}
    \begin{aligned}
    c^{k+1}_1=&\begin{cases}
    \beta c^k_1&\hbox{if }\|A_1y^k+Mz^k+A_2w_0-d\|_\infty>\gamma \big\|A_1y^{k-1}+Mz^{k-1}+A_2w_0-d\|_\infty,\\
    c^k_1&\hbox{otherwise};
\end{cases}\\
c^{k+1}_2=&\begin{cases}
    \beta c^k_2&\hbox{if }\|h^+(z^k,\mu^k,c^k_2)\|_\infty>\gamma \|h^+(z^{k-1},\mu^{k-1},c^{k-1}_2)\|_\infty,\\
    c^k_2&\hbox{otherwise},
\end{cases}
\end{aligned}
\end{equation}}
where $\|\cdot\|_\infty$ denotes the infinity norm.

{Standard convergence results for AL methods apply to convex constrained optimization problems under suitable assumptions on the accuracy with which the inner minimizations are performed; see \cite[Chapter~5]{Ber82a} for a general discussion of convergence and rate of convergence. In practice, the computational performance of the method depends on the choice of the penalty parameters and on the details of the penalty-update rule \eqref{eq:penalty_update}. We examine the sensitivity to these choices computationally in Appendix~\ref{app:algorithmic-tuning}.}

\section{Computational Studies}
\label{sec:experiments}

We have tested the compressed formulation~\eqref{eq:compressed_ta} and the AL method of Section~\ref{sec:alm} on three standard benchmark networks from the Transportation Networks for Research repository~\cite{TNR}: Chicago Sketch, Chicago Regional, and Philadelphia; see Table~\ref{tab:network-characteristics} for selected network statistics. Recall that the AL method consists of \emph{outer iterations}, each of which updates the multipliers and penalty parameters [cf. Eqs.~\eqref{eq:multiplier_update} and \eqref{eq:penalty_update}], and \emph{inner iterations}, which minimize the augmented Lagrangian via L-BFGS-B [cf. Eq.~\eqref{eq:inner_iteration}]; see Section~\ref{sec:alm_iterations}.

\begin{table}[htbp]
\centering
\caption{\color{black}Network characteristics for the original path pools under nominal demand. Singleton paths correspond to OD pairs with a unique path and are excluded from the compressed problem. The expanded path pools considered later are described in Section~\ref{sec:path-pool-expansion}.}
\label{tab:network-characteristics}
{\small
\begin{tabular}{lrrr}
\toprule
 & Chicago Sketch & Chicago Regional & Philadelphia \\
\midrule
Nodes / Links & 933 / 2{,}950 & 12{,}982 / 39{,}018 & 13{,}389 / 40{,}003 \\
Multi-path variables ($n$) & 42{,}774 & 879{,}625 & 1{,}846{,}068 \\
Multi-path OD pairs ($\ell$) & 17{,}464 & 353{,}582 & 601{,}644 \\
Multi-path flow share & 23.5\% & 59.2\% & 70.9\% \\
\bottomrule
\end{tabular}
}
\end{table}

{\color{black}Our computational studies have two distinct purposes. The studies in Sections~\ref{sec:cutoff-thresholds}-\ref{sec:path-pool-expansion} use the nominal demand to study the computational consequences of the compression design parameters, including the threshold $\tau$, rank $r$, and path-pool size. Since the reference minor-path flows are constructed from the nominal solution, the nominal reference solution is representable in the compressed formulation by setting $z=0$. Thus, these experiments primarily examine the dimensional reduction and computational behavior associated with different design choices, while also verifying that the numerical solution remains close to the nominal reference solution.

The study in Section~\ref{sec:demand-perturbations} addresses the situation in which the proposed method is intended to be used. Once the design parameters have been selected from the nominal problem, we keep the path pool, major--minor partition, and SVD basis fixed and apply the compressed formulation to problems with perturbed OD demands. The reference minor-path flows are scaled according to the perturbed OD demands, but neither the partition nor the SVD basis is recomputed. In this scenario, the relevant question is whether the reused major--minor partition and SVD basis continue to provide sufficiently accurate solutions while reducing the computational cost of solving the perturbed problems.}

The two stages also use different initialization protocols, reflecting their different purposes. In Sections~\ref{sec:cutoff-thresholds}-\ref{sec:path-pool-expansion}, both the uncompressed and compressed formulations are initialized from a neutral cold start based on an equal allocation of each OD demand across the available paths. This allows the design parameters to be studied under a common initialization that does not warm-start the optimization variables from the nominal solution. In Section~\ref{sec:demand-perturbations}, by contrast, the nominal path flows are already available from the offline design stage and are therefore reused to initialize the perturbed-demand problems. The resulting experiments represent the intended repeated-use setting rather than a cold-start comparison.

We use three metrics to measure solution accuracy: the BPR gap, defined as the relative difference between the computed and reference objective values; the coefficient of determination $R^2$, which measures agreement between the computed and reference link-flow vectors; and the link travel-time RMSE, which measures the root mean squared error in link travel times.
In particular, let us denote by $x^*$ and $v^*$ the path and link flows of a high-accuracy reference user-equilibrium solution computed on the same path pool and OD demands as the corresponding problem instance. {\color{black}For any computed link-flow vector $\tilde v$, let $\tilde t_i$ and $t_i^*$ denote the travel times on link $i$ corresponding to $\tilde v_i$ and $v_i^*$, respectively, and let
$$
\bar v^*=\frac{1}{m}\sum_{i=1}^m v_i^*.
$$
We use the following three measures of solution accuracy:
\begin{equation}
\label{eq:metric}
\begin{aligned}
\textit{BPR gap}
&=
\frac{f(\tilde v)-f(v^*)}{f(v^*)}\times 100\,(\%),\\
R^2
&=
1-
\frac{\sum_{i=1}^{m}(\tilde v_i-v_i^*)^2}
{\sum_{i=1}^{m}(v_i^*-\bar v^*)^2},\\
\mathrm{RMSE}_t
&=
60\sqrt{\frac{1}{m}\sum_{i=1}^{m}
(\tilde t_i-t_i^*)^2}.
\end{aligned}
\end{equation}
The factor $60$ in the definition of $\mathrm{RMSE}_t$ converts travel times from minutes to seconds.

For the compressed formulation, the computed link-flow vector is
$$
\tilde v
=
B_1'y^*+Dz^*+B_2'w_0+v_0,
$$
where $(y^*,z^*)$ is the solution computed by the AL method. The same three measures are used to evaluate the corresponding uncompressed AL solutions.
}

The main findings from our computational studies are as follows:
{\color{black}\begin{itemize}
\item \emph{Compression yields substantial dimensional reduction.}
Across the three networks, the compressed AL solutions remain in very close agreement with the nominal reference solutions.

\item \emph{Nominal results are insensitive to rank.} Across the tested ranks, the numerical results change very little.

\item \emph{Compression remains effective for expanded path pools.} Large dimensional reductions are maintained for the expanded path pools.

\item \emph{The nominally constructed compression can be reused.} Under perturbed OD demands, the major--minor partition and SVD basis are retained while the reference flows are scaled with demand; high accuracy and substantial computational savings are maintained.

\item \emph{The compressed formulation offers an accuracy-computation tradeoff.} Fixing the minor-path flows at their reference values already yields high accuracy, while allowing low-dimensional adjustments can further improve link travel-time accuracy in more demanding instances.
\end{itemize}}

In what follows, we provide further details on our computational studies. The codes can be found at \url{https://github.com/asu-trans-ai-lab/CompressedTAP}.

\subsection{Experimental Setup}
\label{sec:setup}

For each network, the flow on the path of every singleton OD pair (i.e., an OD pair with a unique path) is fixed to the corresponding OD demand. These paths contribute a constant offset $v_0$ to the link flows; only multi-path OD pairs enter the compressed problem. The compressed problem is solved with the AL iterations of Section~\ref{sec:alm_iterations}, using L-BFGS-B for the inner minimization over $(y, z)$ with $y \geq 0$. The gradient of the AL function with respect to $z$ is computed by the memory-efficient and time-efficient scheme detailed in Appendix~\ref{app:gradient-overhead}, which avoids forming any dense matrix of size $m \times (n-s)$ or $\ell \times (n-s)$. 

Unless stated otherwise, we use $\beta = 10$ and $\gamma = 0.25$, as suggested in \cite[Section~2.2.5]{Ber82a}, initial penalty $c^0 = (10^3, 10^3)$, convergence tolerance $10^{-4}$ on the infinity norm of the constraint violation, and at most 20 outer and 200 inner iterations. The sensitivity of these choices is examined in Appendix~\ref{app:algorithmic-tuning}. The default compression rank {\color{black}$r$ is set to be $20$, with other choices examined in Section~\ref{sec:rank}.}

Solution quality is measured against a high-accuracy reference user-equilibrium solution computed on the same path pool and OD demands as the corresponding problem instance. Details of the reference-solution calculations are given in Appendix~\ref{app:reference-solution}. Accuracy is evaluated using the BPR gap, the coefficient of determination $R^2$, and the link travel-time RMSE defined in Eq.~\eqref{eq:metric}. Computational cost is reported in terms of CPU time and the numbers of outer and inner iterations. For comparison, the uncompressed problem~\eqref{eq:original_ta} is solved with identical algorithmic settings on each network. All CPU times are single-thread measurements on an Apple M1 Pro workstation.

{\color{black}The SVD used to construct the compressed representation is treated as a one-time preprocessing step and is not included in the reported CPU times, which measure only the subsequent optimization solves.}

\subsection{Selection of the Compression Threshold}
\label{sec:cutoff-thresholds}

The threshold $\tau$ controls the partition into $s$ major and $n-s$ minor paths and, together with the rank $r$, determines the dimension of the compressed problem: $s+r$ decision variables replace the original $n$ path-flow variables. We choose $\tau$ at 10 equally spaced quantile levels of the nominal path-flow distribution, so that increasing $\tau$ progressively classifies a larger fraction of low-flow paths as minor. For each threshold, we report the variable reduction $(n-s-r)/n$, the CPU time, the number of outer and inner iterations, the BPR gap, $R^2$, and $\mathrm{RMSE}_t$. Complete results are given in Appendix~\ref{app:full-results}; the distribution of flow across major and minor paths at each threshold is reported in Appendix~\ref{app:pathflow}; the main findings are summarized in Table~\ref{tab:cutoff-summary}.

\begin{figure}[htbp]
    \centering
    \includegraphics[width=0.92\textwidth]{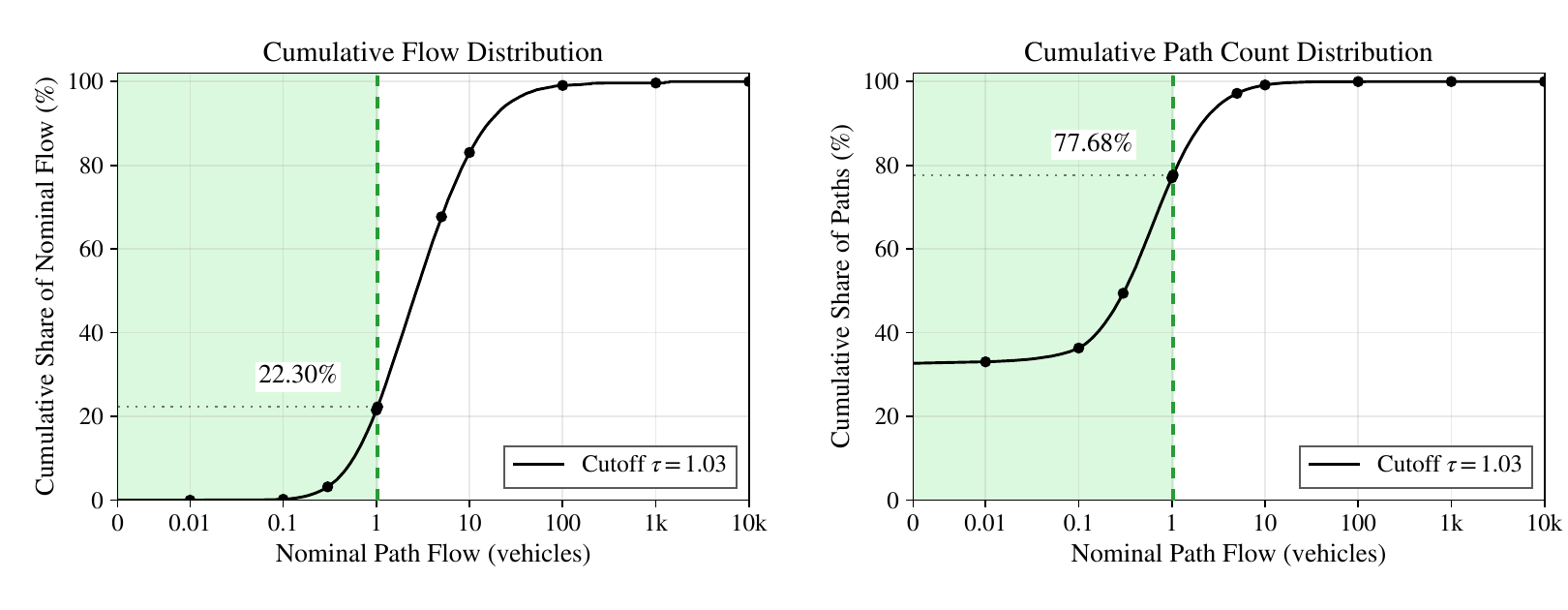}
    \caption{Cumulative distributions of nominal path flow and path count for
    Chicago Regional. The vertical dashed line indicates the cutoff
    $\tau=1.03$. Approximately $77.68\%$ of the paths have nominal flows below
    the cutoff, accounting for $22.30\%$ of the nominal flow associated with
    the multi-path OD pairs. The actual minor-path set is smaller because the
    largest-flow path of each OD pair is retained as a major path regardless
    of the cutoff.}
    \label{fig:path-flow-distribution}
\end{figure}

\begin{figure}[htbp]
    \centering
    \includegraphics[width=0.99\textwidth]{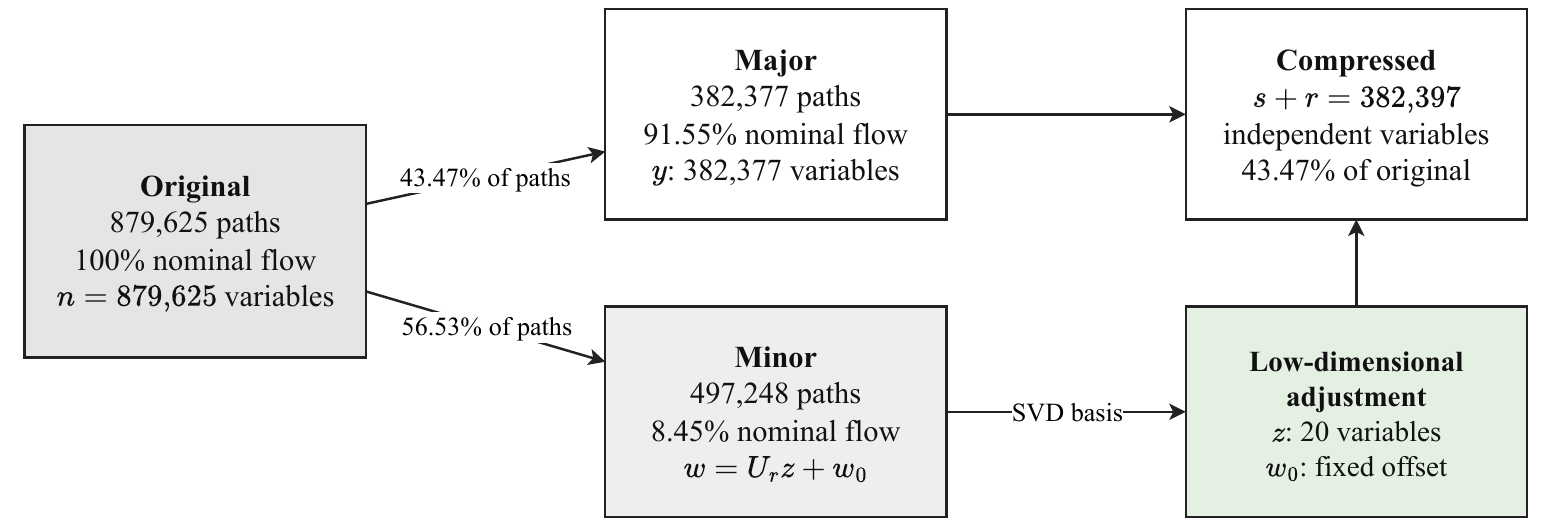}
    \caption{Illustration of the affine compression for Chicago Regional at
    $\tau=1.03$ and $r=20$. The major paths are retained explicitly, while
    the minor-path flows are represented as $w=w_0+U_rz$, where $w_0$ is the
    fixed reference flow and $z$ has dimension $20$. The percentages
    $91.55\%$ and $8.45\%$ denote shares of the nominal flow associated with
    the multi-path OD pairs. The resulting formulation has $382{,}397$
    independent path-flow variables, or $43.47\%$ of the original
    $879{,}625$.}
    \label{fig:affine-compression-example}
\end{figure}

{\color{black}Figure~\ref{fig:path-flow-distribution} illustrates the concentration of nominal path flows for Chicago Regional. At the cutoff $\tau=1.03$, approximately $77.68\%$ of the paths have nominal flows below the cutoff, but these paths account for only $22.30\%$ of the nominal flow associated with the multi-path OD pairs. Thus, a large fraction of the path set consists of relatively low-flow paths. Because our classification rule retains at least one major path for every OD pair, however, not every path below the cutoff is classified as minor. Consequently, the resulting fraction of minor paths is smaller than the cumulative path share shown in the figure; see Fig.~\ref{fig:affine-compression-example}.}

\begin{table}[htbp]
\centering
\caption{Selected cutoff-threshold results across all three networks ($r = 20$). The first row for each network ($\tau = 0$) is the uncompressed baseline. Full tables are in Appendix~\ref{app:full-results}.}
\label{tab:cutoff-summary}
{\small\color{black}
\resizebox{\textwidth}{!}{%
\begin{tabular}{ll rrrrrr}
\toprule
Network & $\tau$ & \shortstack{Variable\\Reduction (\%)} & \shortstack{BPR \\ Gap (\%)} & Link $R^2$ & \shortstack{RMSE (s)} & \shortstack{Outer / Inner\\Iterations} & \shortstack{CPU (s)}  \\
\midrule
Chicago Sketch   & 0  & --- & 0.179 & 0.99562 & 4.123 & 7\,/\,620 & 3.95  \\
                 & 0.46  & 40.02 & 0.029 & 0.99834 & 3.833 & 6\,/\,720 & 3.16  \\
                 & 1.06  & 46.85 & 0.028 & 0.99843 & 3.544 & 8\,/\,721 & 2.99  \\
                 & 4.54  & 54.72 & 0.027 & 0.99854 & 3.952 & 7\,/\,599 & 2.28  \\
\midrule
Chicago Regional & 0   & --- & 0.102 & 0.99939 & 1.653 & 6\,/\,817 & 111.28  \\
                 & 0.23  & 43.69 & 0.009 & 0.99989 & 1.139 & 7\,/\,454 & 51.43  \\
                 & 0.46  & 51.08 & 0.007 & 0.99992 & 1.098 & 6\,/\,616 & 61.77  \\
                 & 1.03  & 56.53 & 0.007 & 0.99993 & 0.986 & 6\,/\,691 & 65.94  \\
\midrule
Philadelphia     & 0  & --- & 0.414 & 0.99559 & 4.717 & 13\,/\,1{,}034 & 315.55  \\
                 & 0.67  & 51.27 & 0.061 & 0.99915 & 3.252 & 8\,/\,851 & 177.08  \\
                 & 0.99  & 54.44 & 0.061 & 0.99916 & 3.252 & 8\,/\,842 & 168.69  \\
                 & 5.33  & 63.61 & 0.054 & 0.99935 & 3.101 & 10\,/\,718 & 137.02  \\
\bottomrule
\end{tabular}%
}}
\end{table}

{\color{black}Table~\ref{tab:cutoff-summary} shows that substantial variable reduction can be achieved while maintaining very high numerical agreement with the nominal reference solution. Across all selected compressed cases, the BPR gap remains below $0.07\%$ and the link-flow $R^2$ exceeds $0.9983$. At the largest selected thresholds, the variable reductions reach $54.7\%$, $56.5\%$, and $63.6\%$ on Chicago Sketch, Chicago Regional, and Philadelphia, respectively, while the corresponding link-flow $R^2$ values are $0.99854$, $0.99993$, and $0.99935$.

Interestingly, the compressed AL solutions are numerically closer to the reference solution than the uncompressed AL solutions in these nominal experiments. For Chicago Sketch, the BPR gap decreases from $0.179\%$ in the uncompressed case to $0.027\%$ at the largest selected threshold; for Chicago Regional, it decreases from $0.102\%$ to $0.007\%$; and for Philadelphia, from $0.414\%$ to $0.054\%$. The corresponding link-flow $R^2$ values also increase. Since the nominal reference solution is representable in the compressed formulation, these differences should be viewed as differences in the numerical behavior of the AL solutions rather than as improvements in the underlying optimal solution produced by compression.

The computational results show a similarly favorable pattern. Every selected compressed case requires less total CPU time than its corresponding uncompressed baseline. At the largest selected thresholds, CPU time decreases from {\color{black}$3.95$~s to $2.28$~s} on Chicago Sketch, from {\color{black}$111.28$~s to $65.94$~s} on Chicago Regional, and from {\color{black}$315.55$~s to $137.02$~s} on Philadelphia, corresponding to reductions of approximately {\color{black}$42\%$, $41\%$, and $57\%$}, respectively. The CPU time is not strictly monotone in $\tau$, because changes in the compressed dimension are accompanied by changes in the numbers of outer and inner iterations, but the reduction in computational cost is substantial throughout the selected cases.

These nominal experiments therefore provide a practical view of how the choice of $\tau$ affects problem dimension and computational behavior. Over the tested range, increasingly aggressive compression yields substantial reductions in dimension while the AL solutions remain in very close agreement with the nominal reference solution. We use these results to identify compression levels that provide substantial dimensional reduction without introducing numerical difficulties; the effectiveness of the resulting representations under changes in demand is examined in Section~\ref{sec:demand-perturbations}.}

{\color{black}Figure~\ref{fig:affine-compression-example} illustrates the resulting compression for Chicago Regional at $\tau=1.03$ and $r=20$. Of the $879{,}625$ path-flow variables, $382{,}377$ correspond to major paths and remain explicit decision variables, while $497{,}248$ correspond to minor paths. The major paths carry $91.55\%$ of the nominal flow associated with the multi-path OD pairs, whereas the minor paths carry $8.45\%$. Rather than optimizing all $497{,}248$ minor-path flows independently, the compressed formulation uses the fixed reference flow $w_0$ together with the $20$-dimensional adjustment $U_r z$. The number of independent path-flow variables is therefore reduced from $879{,}625$ to $s+r=382{,}377+20=382{,}397,$ or $43.47\%$ of the original number.}

\subsection{Selection of the Compression Rank}
\label{sec:rank}

The rank $r$ controls the dimension of the low-dimensional adjustment $U_rz$ to the reference minor-path flows. To examine its effect, we fix the threshold at the 90th percentile of the nominal path-flow distribution on each network and consider $r\in\{5,10,20,50\}$. We report the results for the two larger networks, Chicago Regional and Philadelphia.

\begin{table}[htbp]
\centering
\caption{Effect of rank $r$ on solution quality and computational cost. The uncompressed baseline CPU times are {\color{black}111.28s} (Chicago Regional) and {\color{black}315.55s} (Philadelphia).}
\label{tab:impact-rank}
{\small\color{black}
\begin{tabular}{ll rrrrr}
\toprule
Network & $r$ & \shortstack{BPR \\ Gap (\%)} & Link $R^2$ & \shortstack{RMSE (s)} & \shortstack{Outer / Inner\\Iterations} & \shortstack{CPU (s)}  \\
\midrule
Chicago Regional & 5   & 0.007 & 0.99993 & 0.991 & 9\,/\,491 & 43.66  \\
($\tau = 1.03$, \   & 10  & 0.007 & 0.99992 & 1.023 & 7\,/\,559 & 52.19  \\
56.53\% red.)     & 20  & 0.007 & 0.99993 & 0.986 & 6\,/\,691 & 65.94  \\
                 & 50  & 0.011 & 0.99992 & 1.047 & 9\,/\,1{,}256 & 105.01  \\
\midrule
Philadelphia     & 5   & 0.041 & 0.99946 & 2.630 & 10\,/\,503 & 85.95  \\
($\tau = 5.33$, \   & 10  & 0.049 & 0.99939 & 2.916 & 9\,/\,657 & 120.68  \\
63.61\% red.)     & 20  & 0.054 & 0.99935 & 3.101 & 10\,/\,718 & 137.02  \\
                 & 50  & 0.061 & 0.99929 & 3.322 & 9\,/\,777 & 129.57  \\
\bottomrule
\end{tabular}
}
\end{table}

\begin{figure}[htbp]
    \centering
    \includegraphics[width=\textwidth]{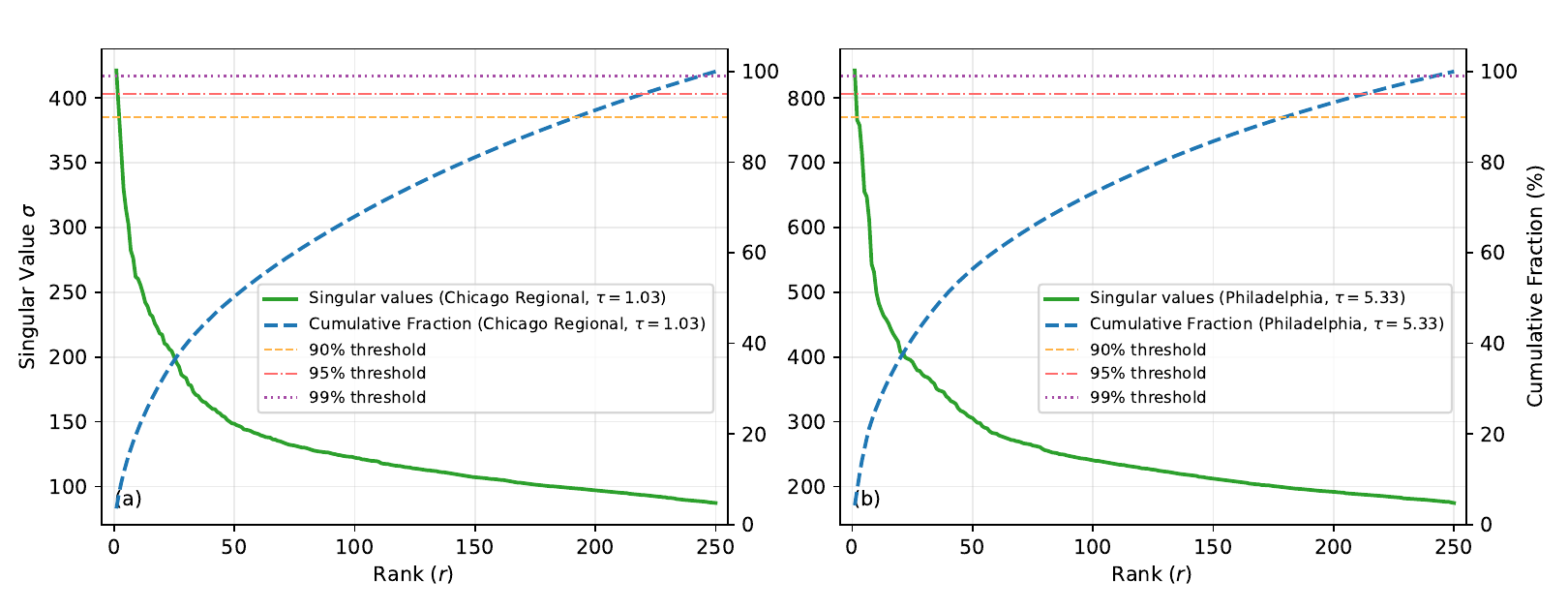}
    \caption{{\color{black}Singular value spectra and cumulative fractions of squared singular values for $B_2$ on Chicago Regional at $\tau = 1.03$ (left) and Philadelphia at $\tau = 5.33$ (right).}}
    \label{fig:svd-spectrum}
\end{figure}

{\color{black}Table~\ref{tab:impact-rank} shows that the numerical solutions for the nominal-demand problems are largely insensitive to the compression rank over the range considered. On Chicago Regional, the BPR gap varies only from $0.007\%$ to $0.011\%$, while the link-flow $R^2$ remains approximately $0.9999$. On Philadelphia, the BPR gap ranges from $0.041\%$ to $0.061\%$, and $R^2$ remains above $0.9992$. The corresponding link travel-time errors also vary only moderately across the tested ranks.

This insensitivity is consistent with the structure of the compressed formulation. For the nominal-demand problem, the reference minor-path flow is the nominal minor-path flow itself, and hence the nominal reference solution is representable with $z=0$ for every value of $r$. The observed differences across ranks therefore mainly reflect the numerical behavior of the AL solution rather than an approximation error caused by the choice of rank.

Figure~\ref{fig:svd-spectrum} provides a complementary perspective. It shows the singular values of $B_2$ together with the cumulative fraction of squared singular values. The gradual spectral decay indicates that a substantially larger rank would be required for high-fidelity reconstruction of $B_2$ itself. The nominal computations, however, remain stable over a much smaller range of ranks. Thus, spectral reconstruction accuracy alone is not an appropriate criterion for selecting $r$; the practical value of the low-dimensional adjustment is better assessed when the compressed representation is reused for problems that differ from the nominal one.

The computational results provide an additional consideration in selecting $r$. Although CPU time is not monotone in $r$, the largest rank tested gives the largest BPR gap and travel-time RMSE on both networks and the largest CPU time on Chicago Regional. All ranks up to $r=20$ remain substantially faster than the corresponding uncompressed baselines.

We use $r=20$ as a moderate intermediate value in the remaining experiments. This choice is not intended to be optimal; rather, it provides a representative rank within the range over which the nominal numerical results are stable.}

\subsection{\color{black}Effect of Path-Pool Expansion}
\label{sec:path-pool-expansion}

{\color{black}The {\color{black}expansion} of the path pool controls the number of candidate routes represented for each OD pair and, consequently, the dimension of the path-based formulation. To examine whether the compression framework remains effective for expanded path pools, we construct two nested expansions, denoted by X4 and X8. The original candidate-path pool is retained as a strict subset, and additional $k$-shortest-path candidates are added only for the original multi-path OD pairs until each such pair has at least four or eight represented paths, respectively. OD pairs having a single original path remain fixed through $v_0$.

For each expanded pool, we compute a new nominal reference assignment, which provides the nominal path flows used to construct the major--minor partition and the corresponding reference link flows. We use the smallest positive cutoff selected for each network in Table~\ref{tab:cutoff-summary}: $\tau=0.461887$ for Chicago Sketch, $0.226544$ for Chicago Regional, and $0.666284$ for Philadelphia. We fix $r=20$ and compare the uncompressed and compressed formulations using the same cold-start AL settings as in the cutoff-threshold experiments. The SVD is computed once for each expanded pool and is not included in the reported optimization times.}

\begin{table}[htbp]
\centering
\caption{{\color{black}Effect of path-pool {\color{black}expansion} on solution quality and computational cost ($r=20$).}}
\label{tab:path-pool-expansion}
{\scriptsize\color{black}
\resizebox{\textwidth}{!}{%
\begin{tabular}{lrrrrrr}
\toprule
Network & \shortstack{Variable\\Reduction (\%)} & \shortstack{BPR Gap (\%)\\Uncomp. / Comp.} & \shortstack{Link $R^2$\\Uncomp. / Comp.} & \shortstack{RMSE (s)\\Uncomp. / Comp.} & \shortstack{CPU (s)\\Uncomp. / Comp.} & \shortstack{CPU\\Reduction (\%)}  \\
\midrule
Chicago Sketch X4  & 63.81 & $0.887 / 0.023$ & $0.99363 / 0.99868$ & $5.235 / 3.818$ & $8.59 / 3.49$ & 59.41  \\
Chicago Sketch X8  & 81.67 & $1.700 / 0.025$ & $0.99253 / 0.99859$ & $5.996 / 4.218$ & $21.47 / 6.63$ & 69.13  \\
\midrule
Chicago Regional X4  & 65.73 & $0.315 / 0.008$ & $0.99905 / 0.99991$ & $2.478 / 0.981$ & $239.36 / 87.73$ & 63.35  \\
Chicago Regional X8  & 82.63 & $0.369 / 0.009$ & $0.99905 / 0.99991$ & $1.976 / 1.092$ & $573.35 / 137.91$ & 75.95  \\
\midrule
Philadelphia X4  & 65.38 & $1.121 / 0.060$ & $0.98520 / 0.99910$ & $7.248 / 3.029$ & $569.95 / 180.29$ & 68.37  \\
Philadelphia X8  & 81.40 & $2.320 / 0.059$ & $0.96530 / 0.99915$ & $7.344 / 3.449$ & $1148.75 / 251.18$ & 78.13  \\
\bottomrule
\end{tabular}%
}}
\end{table}

{\color{black}Table~\ref{tab:path-pool-expansion} shows that substantial dimensional reduction is maintained as the path pools are expanded. The variable reduction ranges from approximately $64\%$ to $83\%$, while the compressed AL solutions remain in very close agreement with the corresponding nominal reference solutions. Across all six cases, the compressed BPR gap is at most $0.060\%$, and the compressed link-flow $R^2$ remains above $0.9985$. Moreover, the numerical agreement is stable from X4 to X8 for each network.

Under the common cold-start stopping rule, the compressed AL solutions are also substantially closer to the reference solutions than the corresponding uncompressed AL endpoints, particularly for the larger path pools. For example, on Philadelphia X8, the uncompressed BPR gap is $2.320\%$ and the link-flow $R^2$ is $0.96530$, whereas the compressed formulation gives a BPR gap of $0.059\%$ and an $R^2$ of $0.99915$. Similar patterns occur on the two Chicago networks. Since the nominal reference solution is representable in the compressed formulation, these differences should be interpreted as differences in the numerical behavior of the two AL solves under the common stopping protocol, rather than as an intrinsic improvement in the underlying optimum due to compression.

The computational advantage also increases with path-pool {\color{black}expansion}. The compressed formulation reduces CPU time by {\color{black}$59.4\%$--$68.4\%$} for the X4 pools and by {\color{black}$69.1\%$--$78.1\%$} for the X8 pools. The larger reductions for X8 are consistent with the fact that the compressed formulation becomes increasingly advantageous as the number of path variables grows. However, because the uncompressed and compressed runs do not always terminate at comparable levels of solution accuracy, these CPU reductions should be interpreted as comparisons under a common stopping protocol rather than as equal-accuracy computational comparisons.

Overall, these results show that the computational benefit of the compression framework is not tied to a particular path-pool size. As the number of candidate paths increases, the compressed formulation achieves increasingly large reductions in dimension and computational time, while the computed solutions remain in close numerical agreement with the nominal reference solutions.}

\subsection{\color{black}Performance under Demand Perturbations}
\label{sec:demand-perturbations}

{\color{black}We next evaluate the setting in which the proposed compression framework is intended to be used. Starting from the nominal problem, we keep the path pool, major--minor partition, and rank-$20$ SVD basis fixed, and perturb the multi-path OD demands by up to $\pm5\%$ and $\pm10\%$. For each perturbed instance, the reference minor-path flows are scaled by the same factors as the corresponding OD demands, while the major--minor partition and SVD basis are not recomputed.

In this repeated-use setting, the nominal path flows are naturally available. For each perturbed OD pair, we scale its nominal path flows by the same factor as its demand. These scaled flows are used to initialize the uncompressed formulation. For the compressed formulation, the scaled minor-path flows define the reference flow $w_0$, the retained major-path variables are initialized at their correspondingly scaled nominal values, and $z$ is initialized at zero. Thus, both formulations reuse only information available from the nominal problem.

For the original path pools, we use the largest selected thresholds in Table~\ref{tab:cutoff-summary}. For the X4 and X8 pools, we retain the network-specific partitions constructed in Section~\ref{sec:path-pool-expansion}. Following the setup of Section~\ref{sec:setup}, solution quality is evaluated relative to the corresponding reference solution for each perturbed-demand instance.

Table~\ref{tab:demand-perturbations-original} shows that the reused compressed formulation remains highly accurate on the original path pools. Across all six network--perturbation cases, the compressed BPR gap is at most {\color{black}$0.002\%$}, while the compressed link-flow $R^2$ is at least {\color{black}$0.99997$}. The difference from the corresponding uncompressed solutions is therefore very small, despite variable reductions of approximately $55\%$--$64\%$. The results are also stable as the perturbation magnitude increases from $\pm5\%$ to $\pm10\%$, indicating that the compressed formulation remains effective when the nominally constructed partition and SVD basis are reused over the tested range of demand changes. 

The reused compressed formulation also provides a consistent computational advantage. On the original path pools, the compressed formulation requires less CPU time in all six cases, with reductions ranging from {\color{black}$58.5\%$ to $77.2\%$}. Thus, reuse of the nominally constructed partition and SVD basis preserves essentially the same solution quality while substantially reducing the cost of solving the perturbed problems.

}

\begin{table}[htbp]
\centering
\caption{{\color{black}Demand perturbation results on the original path pools ($r=20$). Entries are arithmetic means across demand seeds 10, 28, and 42.}}
\label{tab:demand-perturbations-original}
{\scriptsize\color{black}
\resizebox{\textwidth}{!}{%
\begin{tabular}{lrrrrrrrr}
\toprule
Network & Pert. & \shortstack{Variable\\Reduction (\%)} & \shortstack{BPR Gap (\%)\\Uncomp. / Comp.} & \shortstack{Link $R^2$\\Uncomp. / Comp.} & \shortstack{RMSE (s)\\Uncomp. / Comp.} & \shortstack{CPU (s)\\Uncomp. / Comp.} & \shortstack{CPU\\Reduction (\%)}  \\
\midrule
Chicago Sketch   & $\pm5\%$   & 54.72 & $0.001 / 0.001$ & $0.99998 / 0.99998$ & $0.315 / 0.328$ & $2.39 / 0.77$ & 67.75  \\
                 & $\pm10\%$  & 54.72 & $0.001 / 0.001$ & $0.99998 / 0.99997$ & $0.393 / 0.537$ & $2.68 / 0.79$ & 70.39  \\
Chicago Regional & $\pm5\%$   & 56.53 & $0.002 / 0.002$ & $0.99998 / 0.99998$ & $0.307 / 0.440$ & $83.88 / 34.83$ & 58.47  \\
                 & $\pm10\%$  & 56.53 & $0.002 / 0.002$ & $0.99998 / 0.99999$ & $0.320 / 0.432$ & $85.26 / 29.07$ & 65.90  \\
Philadelphia     & $\pm5\%$   & 63.61 & $0.002 / 0.002$ & $0.99999 / 0.99999$ & $0.282 / 0.323$ & $181.32 / 41.38$ & 77.18  \\
                 & $\pm10\%$  & 63.61 & $0.002 / 0.002$ & $0.99998 / 0.99998$ & $0.453 / 0.483$ & $160.74 / 57.07$ & 64.50  \\
\bottomrule
\end{tabular}%
}}
\end{table}

{\color{black}We next apply the same reuse procedure to the {\color{black}expanded} X4 and X8 path pools. Table~\ref{tab:demand-perturbations-expansion} summarizes the results, with complete seed-level results reported in Appendix~\ref{app:demand-perturbation-full}. The accuracy is particularly stable in these experiments. Across all twelve cases, the compressed BPR gap is at most {\color{black}$0.006\%$}, and the compressed link-flow $R^2$ is at least {\color{black}$0.99990$}, even though the variable reduction ranges from approximately $64\%$ to $83\%$. The compressed and uncompressed solutions are therefore virtually indistinguishable in terms of both objective value and link flows.}

\begin{table}[htbp]
\centering
\caption{\color{black}Demand perturbation results on the X4 and X8 pools ($r=20$). Entries are arithmetic means across demand seeds 10, 28, and 42.}
\label{tab:demand-perturbations-expansion}
{\scriptsize\color{black}
\resizebox{\textwidth}{!}{%
\begin{tabular}{lrrrrrrr}
\toprule
Network & Pert. & \shortstack{Variable\\Reduction (\%)} & \shortstack{BPR Gap (\%)\\Uncomp. / Comp.} & \shortstack{Link $R^2$\\Uncomp. / Comp.} & \shortstack{RMSE (s)\\Uncomp. / Comp.} & \shortstack{CPU (s)\\Uncomp. / Comp.} & \shortstack{CPU\\Reduction (\%)}  \\
\midrule
Chicago Sketch X4 & $\pm5\%$   & 63.81 & $0.002 / 0.002$ & $0.99997 / 0.99997$ & $0.363 / 0.395$ & $3.34 / 0.93$ & 72.15  \\
                  & $\pm10\%$  & 63.81 & $0.002 / 0.002$ & $0.99997 / 0.99997$ & $0.444 / 0.681$ & $3.37 / 1.12$ & 66.66  \\
Chicago Sketch X8 & $\pm5\%$   & 81.67 & $0.002 / 0.002$ & $0.99997 / 0.99997$ & $0.380 / 0.677$ & $5.72 / 1.39$ & 75.75  \\
                  & $\pm10\%$  & 81.67 & $0.002 / 0.002$ & $0.99996 / 0.99996$ & $0.470 / 0.639$ & $5.74 / 1.56$ & 72.81  \\
Chicago Regional X4 & $\pm5\%$   & 65.73 & $0.003 / 0.003$ & $0.99998 / 0.99998$ & $0.343 / 0.382$ & $155.47 / 71.51$ & 54.00  \\
                    & $\pm10\%$  & 65.73 & $0.003 / 0.003$ & $0.99998 / 0.99998$ & $0.353 / 0.352$ & $154.24 / 71.48$ & 53.66  \\
Chicago Regional X8 & $\pm5\%$   & 82.63 & $0.004 / 0.003$ & $0.99998 / 0.99998$ & $0.352 / 0.361$ & $280.36 / 136.58$ & 51.28  \\
                    & $\pm10\%$  & 82.63 & $0.004 / 0.003$ & $0.99997 / 0.99998$ & $0.367 / 0.390$ & $280.06 / 133.70$ & 52.26  \\
Philadelphia X4 & $\pm5\%$   & 65.38 & $0.004 / 0.004$ & $0.99994 / 0.99994$ & $0.383 / 0.479$ & $251.10 / 109.99$ & 56.20  \\
                & $\pm10\%$  & 65.38 & $0.005 / 0.005$ & $0.99993 / 0.99993$ & $0.502 / 0.616$ & $253.54 / 101.89$ & 59.81  \\
Philadelphia X8 & $\pm5\%$   & 81.40 & $0.005 / 0.005$ & $0.99992 / 0.99992$ & $0.394 / 0.602$ & $473.32 / 174.50$ & 63.13  \\
                & $\pm10\%$  & 81.40 & $0.005 / 0.006$ & $0.99991 / 0.99991$ & $0.536 / 0.716$ & $487.15 / 178.56$ & 63.35  \\
\bottomrule
\end{tabular}%
}}
\end{table}

{\color{black}The computational reductions remain substantial for the expanded path pools. For X4, CPU time is reduced by approximately 53.7\%–72.2\%, while for X8 the reduction ranges from approximately 51.3\% to 75.8\%. The savings vary across networks and do not increase monotonically with path-pool size.

Note that since these experiments reuse nominal-flow information for initialization, whereas Sections~\ref{sec:cutoff-thresholds}-\ref{sec:path-pool-expansion} use neutral cold starts, their absolute accuracy, iteration counts, and CPU times should not be compared directly across the two stages.

Overall, these results support the intended repeated-use interpretation of the proposed framework. A compressed representation constructed once from nominal path-flow information can be reused for perturbed OD demands without recomputing either the major--minor partition or the SVD. Across both the original and expanded path pools, the reused major--minor partition and SVD basis continue to provide extremely high solution accuracy while substantially reducing computational cost.}

{\color{black}We further compare the compressed formulation with two reduced alternatives under the same demand perturbations. Major-only sets all minor-path flows to zero and optimizes only the major-path flows. Fixed sets the minor-path flows equal to the reference flow $w_0$ and optimizes only the major-path flows. Compressed uses
$$ w=U_rz+w_0 $$
and therefore allows low-dimensional adjustments to the reference minor-path flows. All three methods use the same major--minor partition and warm-start information. Table~\ref{tab:perturbed-three-arms} reports arithmetic means over seeds 10, 28, and 42.}

\begin{table}[htbp]
\centering
\caption{\color{black}Major-only, Fixed, and Compressed on the original pools under demand perturbations ($r=20$). Triplets follow this method order.}
\label{tab:perturbed-three-arms}
{\scriptsize\color{black}
\resizebox{\textwidth}{!}{%
\begin{tabular}{llrrrr}
\toprule
Network & Pert. & BPR Gap (\%) & Link $R^2$ & \shortstack{RMSE (s)} & CPU (s) \\
\midrule
Chicago Sketch & $\pm5\%$ & 0.004 / 0.001 / 0.001 & 0.99988 / 0.99998 / 0.99998 & 1.972 / 0.333 / 0.328 & 0.33 / 0.10 / 0.77 \\
Chicago Sketch & $\pm10\%$ & 0.004 / 0.001 / 0.001 & 0.99985 / 0.99997 / 0.99997 & 2.012 / 0.548 / 0.537 & 0.26 / 0.09 / 0.79 \\
Chicago Regional & $\pm5\%$ & 0.005 / 0.002 / 0.002 & 0.99993 / 0.99999 / 0.99998 & 0.978 / 0.548 / 0.440 & 12.94 / 7.08 / 34.83 \\
Chicago Regional & $\pm10\%$ & 0.005 / 0.002 / 0.002 & 0.99993 / 0.99998 / 0.99999 & 1.000 / 0.541 / 0.432 & 11.01 / 10.51 / 29.07 \\
Philadelphia & $\pm5\%$ & 0.004 / 0.002 / 0.002 & 0.99993 / 0.99999 / 0.99999 & 0.841 / 0.327 / 0.323 & 22.23 / 9.52 / 41.38 \\
Philadelphia & $\pm10\%$ & 0.005 / 0.003 / 0.002 & 0.99993 / 0.99998 / 0.99998 & 0.882 / 0.627 / 0.483 & 26.40 / 9.36 / 57.07 \\
\bottomrule
\end{tabular}%
}}
\end{table}

{\color{black}The comparison illustrates the accuracy--computation tradeoff provided by the compressed formulation. Major-only is consistently less accurate than the methods that retain the reference minor-path flows. Fixed already provides high accuracy at very low computational cost. Allowing the additional low-dimensional adjustment further reduces the mean travel-time RMSE in all six cases. The improvement over Fixed is small on Chicago Sketch and on Philadelphia at $\pm5\%$, but is approximately $20\%$ on Chicago Regional and $23\%$ on Philadelphia at $\pm10\%$. These additional accuracy gains come at higher computational cost.}

\section{Conclusion and Discussion}
\label{sec:conclusion}

{\color{black}We developed a compressed path-based formulation for large-scale traffic assignment. Major paths are retained explicitly, while reference flows are used for the minor paths and deviations from these flows are restricted to a low-dimensional subspace obtained from a truncated singular value decomposition. The resulting formulation preserves convexity, and feasibility is ensured by appropriately scaling the reference minor-path flows when the OD demands change. We solve the compressed problem using an AL method with separate penalty parameters for different constraint groups and adaptive penalty updates.

Our computational studies first used nominal-demand problems to examine the design and computational behavior of the compressed formulation. Across the Chicago Sketch, Chicago Regional, and Philadelphia networks, substantial dimensionality reduction is achieved while the computed solutions remain in close agreement with the nominal reference solutions. The results are largely insensitive to the compression rank over the range tested, despite the gradual decay of the singular values of the minor path--link incidence matrix. The dimensional and computational advantages of the compression are also maintained for expanded path pools. These results show that spectral reconstruction accuracy alone is not an appropriate criterion for selecting the compression parameters.

We then considered the intended repeated-use setting. The major--minor partition and SVD basis constructed from the nominal problem are reused while the OD demands are perturbed, with the reference minor-path flows scaled according to the new demands. The compressed formulation remains highly accurate without recomputing either the partition or the SVD basis. Computational time is reduced by approximately $58\%$--$77\%$ on the original path pools and $51\%$--$76\%$ on the expanded path pools. Comparisons with simpler reduced formulations also illustrate an accuracy--computation tradeoff. Fixing the minor-path flows at their reference values already provides high accuracy, while allowing low-dimensional adjustments can further improve link travel-time accuracy in more demanding instances.

More broadly, the framework developed here may apply to other large-scale resource allocation problems of the form
\begin{align*}
    \min_{x,v}\quad & f(v)\\
    \text{s.t.}\quad &v=B'x+v_0,\ \ Ax=d, \ \ x\geq 0,
\end{align*}
where $x$ represents high-dimensional decision variables, $v$ represents lower-dimensional aggregate variables, and $B$ describes the mapping from the high-dimensional decision variables to the aggregate variables. Our results suggest that useful compression need not require strong low-rank structure of the full representation. Instead, nominal solution information can be used to identify dominant variables that should be retained explicitly and to construct reference values for the remaining variables, whose deviations can then be represented in a lower-dimensional space. Exploring this idea in other network and resource allocation problems is a promising direction for future work.}

\bibliographystyle{alpha}
\bibliography{ref}

@book{bertsekas2016nonlinear,
  title={Nonlinear Programming},
  author={Bertsekas, Dimitri P.},
  year={2016},
  edition = {3},
  publisher={Athena Scientific}
}

@book{beckmann1956studies,
  title={Studies in the economics of transportation},
  author={Beckmann, Martin and McGuire, C. and Winsten, C.},
  year={1956},
  publisher={Yale University Press}
}

@book{Ber82a,
  author = {Bertsekas, Dimitri. P.},
  title = {Constrained Optimization and Lagrange Multiplier Methods},
  publisher = {Academic Press},
  year = {1982},
  note = {Republished by Athena Scientific, Belmont, MA, 1997}
}

@article{frank1956algorithm,
  title={An algorithm for quadratic programming},
  author={Frank, Marguerite and Wolfe, Philip},
  journal={Naval Research Logistics Quarterly},
  volume={3},
  number={1--2},
  pages={95--110},
  year={1956},
  publisher={Wiley}
}

@article{jayakrishnan1994fast,
  title={A Faster Path-Based Algorithm for Traffic Assignment},
  author={Jayakrishnan, R. and Tsai, W. K. and Prashker, J. and Rajadhyaksha, S.},
  journal={Transportation Research Record},
  volume={1443},
  pages={75--83},
  year={1994}
}

@techreport{bertsekas1980class,
  title={A Class of Optimal Routing Algorithms for Communication Networks},
  author={Bertsekas, Dimitri P.},
  year={1980},
  institution={Massachusetts Institute of Technology, Laboratory for Information and Decision Systems},
  number={LIDS-P-1042},
  address={Cambridge, MA},
  note={Presented at ICCC 80, Atlanta, GA, Oct. 1980}
}

@techreport{bertsekas2011centralized,
  title={Centralized and Distributed Newton Methods for Network Optimization and Extensions},
  author={Bertsekas, Dimitri P.},
  year={2011},
  month={April},
  institution={Massachusetts Institute of Technology, Laboratory for Information and Decision Systems},
  number={LIDS-2866},
  address={Cambridge, MA}
}

@article{mitradjieva2013stiff,
  title={The stiff is Moving—{Conjugate} direction {Frank-Wolfe} methods with applications to traffic assignment},
  author={Mitradjieva, Maria and Lindberg, Per Olov},
  journal={Transportation Science},
  volume={47},
  number={2},
  pages={280--293},
  year={2013},
  publisher={INFORMS}
}

@book{sheffi1985urban,
  title={Urban Transportation Networks: Equilibrium Analysis with Mathematical Programming Methods},
  author={Sheffi, Yosef},
  year={1985},
  publisher={Prentice-Hall},
  address={Englewood Cliffs, NJ}
}

@book{ahuja1993network,
  title={Network Flows: Theory, Algorithms, and Applications},
  author={Ahuja, Ravindra K. and Magnanti, Thomas L. and Orlin, James B.},
  year={1993},
  publisher={Prentice Hall},
  address={Englewood Cliffs, NJ}
}

@book{patriksson1994traffic,
  title={The Traffic Assignment Problem: Models and Methods},
  author={Patriksson, Michael},
  year={1994},
  publisher={VSP},
  address={Utrecht, The Netherlands}
}

@book{bertsekas1998network,
  title={Network Optimization: Continuous and Discrete Models},
  author={Bertsekas, Dimitri P.},
  year={1998},
  publisher={Athena Scientific},
  address={Belmont, MA}
}

@book{bertsekas2009convex,
  title={Convex Optimization Theory},
  author={Bertsekas, Dimitri P.},
  year={2009},
  publisher={Athena Scientific},
  address={Belmont, MA}
}

@article{gallager1977minimum,
  title={A minimum delay routing algorithm using distributed computation},
  author={Gallager, Robert},
  journal={IEEE Transactions on Communications},
  volume={25},
  number={1},
  pages={73--85},
  year={1977},
  publisher={IEEE}
}

@article{bertsekas1984second,
  title={Second derivative algorithms for minimum delay distributed routing in networks},
  author={Bertsekas, Dimitri and Gafni, Eli and Gallager, Robert},
  journal={IEEE Transactions on Communications},
  volume={32},
  number={8},
  pages={911--919},
  year={1984},
  publisher={IEEE}
}

@book{bertsekas2019reinforcement,
  title={Reinforcement Learning and Optimal Control},
  author={Bertsekas, Dimitri P.},
  year={2019},
  publisher={Athena Scientific Belmont, MA}
}

@book{bertsekas2012dynamic,
  title={Dynamic Programming and Optimal Control},
  author={Bertsekas, Dimitri P.},
  volume={2},
  edition={4},
  year={2012},
  publisher={Athena Scientific}
}

@article{boyles2020transportation,
  title={Transportation network analysis},
  author={Boyles, Stephen D. and Lownes, Nicholas E. and Unnikrishnan, Avinash},
  journal={Volume I: Static and Dynamic Traffic Assignment},
  year={2020}
}

@article{bar2002origin,
  title={Origin-based algorithm for the traffic assignment problem},
  author={Bar-Gera, Hillel},
  journal={Transportation Science},
  volume={36},
  number={4},
  pages={398--417},
  year={2002},
  publisher={INFORMS}
}

@article{nie2010class,
  title={A class of bush-based algorithms for the traffic assignment problem},
  author={Nie, Yu Marco},
  journal={Transportation Research Part B: Methodological},
  volume={44},
  number={1},
  pages={73--89},
  year={2010},
  publisher={Elsevier}
}

@article{xie2018greedy,
  title={A greedy path-based algorithm for traffic assignment},
  author={Xie, Jun and Nie, Yu and Liu, Xiaobo},
  journal={Transportation Research Record},
  volume={2672},
  number={48},
  pages={36--44},
  year={2018},
  publisher={SAGE Publications}
}

@article{liu2023alternating,
  title={An alternating direction method of multipliers for solving user equilibrium problem},
  author={Liu, Zhiyuan and Chen, Xinyuan and Hu, Jintao and Wang, Shuaian and Zhang, Kai and Zhang, Honggang},
  journal={European Journal of Operational Research},
  volume={310},
  number={3},
  pages={1072--1084},
  year={2023},
  publisher={Elsevier}
}

@article{liu2024novel,
  title={A novel parallel computing framework for traffic assignment problem},
  author={Liu, Zhiyuan and Dong, Yu and Zhang, Honggang and Zheng, Nan and Huang, Kai},
  journal={Transportation Research Part E: Logistics and Transportation Review},
  volume={189},
  pages={103687},
  year={2024},
  publisher={Elsevier}
}

@article{liu2025admm,
  title={An {ADMM}-based splitting algorithm with multiplier sequential updates for solving traffic assignment},
  author={Liu, Pengjie and Shao, Hu and Xu, Shengbei and Fainman, Emily Zhu and Tang, Chunkai},
  journal={Transportmetrica A: Transport Science},
  pages={1--33},
  year={2025},
  publisher={Taylor \& Francis}
}

@misc{zhou2014dtalite,
  title={DTALite: A queue-based mesoscopic traffic simulator for fast model evaluation and calibration},
  author={Zhou, Xuesong and Taylor, Jeffrey},
  year={2014},
  publisher={Taylor \& Francis}
}

@article{zhou2022meso,
  title={A meso-to-macro cross-resolution performance approach for connecting polynomial arrival queue model to volume-delay function with inflow demand-to-capacity ratio},
  author={Zhou, Xuesong and Cheng, Qixiu and Wu, Xin and Li, Peiheng and Belezamo, Baloka and Lu, Jiawei and Abbasi, Mohammad},
  journal={Multimodal Transportation},
  volume={1},
  number={2},
  pages={100017},
  year={2022}
}

@misc{Path4GMNS,
  title={{Path4GMNS}},
  author={Li, Peiheng and Zhou, Xuesong},
  journal={GitHub},
  year={2025},
  note={Retrieved from \url{https://github.com/jdlph/Path4GMNS}},
}

@article{camargo2015aequilibrae,
  title={AequilibraE: a free QGIS add-on for transportation modeling},
  author={Camargo, Pedro},
  journal={Foss4g North America},
  year={2015}
}

@article{peeta1999generalized,
  title={Generalized singular value decomposition approach for consistent on-line dynamic traffic assignment},
  author={Peeta, Srinivas and Bulusu, Srinivas},
  journal={Transportation Research Record},
  volume={1667},
  number={1},
  pages={77--87},
  year={1999},
  publisher={SAGE Publications}
}

@article{hu2009identification,
  title={Identification of vehicle sensor locations for link-based network traffic applications},
  author={Hu, Shou-Ren and Peeta, Srinivas and Chu, Chun-Hsiao},
  journal={Transportation Research Part B: Methodological},
  volume={43},
  number={8-9},
  pages={873--894},
  year={2009},
  publisher={Elsevier}
}

@article{zhou2010information,
  title={An information-theoretic sensor location model for traffic origin-destination demand estimation applications},
  author={Zhou, Xuesong and List, George F.},
  journal={Transportation Science},
  volume={44},
  number={2},
  pages={254--273},
  year={2010},
  publisher={INFORMS}
}

@article{zhang2025lowrank,
  title={Low-rank approximation of path-based traffic network models},
  author={Zhang, Pengji and Qian, Sean},
  journal={Transportation Research Part C},
  volume={150},
  pages={104151},
  year={2025}
}

@article{eckart1936approximation,
  title={The approximation of one matrix by another of lower rank},
  author={Eckart, Carl and Young, Gale},
  journal={Psychometrika},
  volume={1},
  number={3},
  pages={211--218},
  year={1936},
  publisher={Springer}
}

@article{golub1965svd,
  title={Calculating the singular values and pseudo-inverse of a matrix},
  author={Golub, Gene H. and Kahan, William},
  journal={Journal of the Society for Industrial and Applied Mathematics, Series B: Numerical Analysis},
  volume={2},
  number={2},
  pages={205--224},
  year={1965},
  publisher={SIAM}
}

@article{sirovich1987low,
  title={Low-dimensional procedure for the characterization of human faces},
  author={Sirovich, Lawrence and Kirby, Michael},
  journal={Journal of the Optical Society of America A},
  volume={4},
  number={3},
  pages={519--524},
  year={1987},
  publisher={Optica}
}

@article{markowitz1952portfolio,
  title={Portfolio selection},
  author={Markowitz, Harry},
  journal={The Journal of Finance},
  volume={7},
  number={1},
  pages={77--91},
  year={1952},
  publisher={JSTOR}
}

@article{ross1976arbitrage,
  title={The arbitrage theory of capital asset pricing},
  author={Ross, Stephen A.},
  journal={Journal of Economic Theory},
  volume={13},
  number={3},
  pages={341--360},
  year={1976},
  publisher={Elsevier}
}

@article{chopra1993effect,
  title={The effect of errors in means, variances, and covariances on optimal portfolio choice},
  author={Chopra, Vijay Kumar and Ziemba, William T.},
  journal={Journal of Portfolio Management},
  volume={19},
  number={2},
  pages={6--11},
  year={1993}
}

@book{vapnik1995nature,
  title={The Nature of Statistical Learning Theory},
  author={Vapnik, Vladimir N.},
  year={1995},
  publisher={Springer}
}

@article{scholkopf1998nonlinear,
  title={Nonlinear component analysis as a kernel eigenvalue problem},
  author={Sch{\"o}lkopf, Bernhard and Smola, Alexander and M{\"u}ller, Klaus-Robert},
  journal={Neural Computation},
  volume={10},
  number={5},
  pages={1299--1319},
  year={1998},
  publisher={MIT Press}
}

@inproceedings{rahimi2008random,
  title={Random features for large-scale kernel machines},
  author={Rahimi, Ali and Recht, Benjamin},
  booktitle={Advances in Neural Information Processing Systems},
  volume={20},
  pages={1177--1184},
  year={2008}
}

@inproceedings{sarwar2000application,
  title={Application of dimensionality reduction in recommender system---a case study},
  author={Sarwar, Badrul and Karypis, George and Konstan, Joseph and Riedl, John},
  booktitle={Proceedings of the ACM WebKDD Workshop},
  pages={82--90},
  year={2000}
}

@article{koren2009matrix,
  title={Matrix factorization techniques for recommender systems},
  author={Koren, Yehuda and Bell, Robert and Volinsky, Chris},
  journal={Computer},
  volume={42},
  number={8},
  pages={30--37},
  year={2009},
  publisher={IEEE}
}

@article{candes2009exact,
  title={Exact matrix completion via convex optimization},
  author={Cand{\`e}s, Emmanuel J. and Recht, Benjamin},
  journal={Foundations of Computational Mathematics},
  volume={9},
  number={6},
  pages={717--772},
  year={2009},
  publisher={Springer}
}

@article{prakash2017reducing,
  title={Reducing the dimension of online calibration in dynamic traffic assignment systems},
  author={Prakash, A. Arun and Seshadri, Ravi and Antoniou, Constantinos and Pereira, Francisco C. and Ben-Akiva, Moshe E.},
  journal={Transportation Research Record},
  volume={2667},
  number={1},
  pages={96--107},
  year={2017},
  publisher={SAGE Publications}
}

@book{golub2013matrix,
  title={Matrix Computations},
  author={Golub, Gene H. and Van Loan, Charles F.},
  edition={4},
  publisher={Johns Hopkins University Press},
  year={2013}
}

@misc{TNR, 
author = {{Transportation Networks for Research Core Team}}, title = {Transportation Networks for Research}, 
note = {Accessed September 2026}
}

\appendix

\section{Notation}\label{app:notation}

\begin{longtable}{ll}
\caption{Summary of Main Notation} \label{tab:notation} \\
\toprule
Symbol & Description \\
\midrule
\endfirsthead
\toprule
Symbol & Description \\
\midrule
\endhead
\midrule
\endfoot
\bottomrule
\endlastfoot
\multicolumn{2}{l}{\textbf{Network and Demand}} \\
$d$ & OD demand vector (multi-path ODs only) \\
$v$ & Link flow vector \\
$v_0$ & Link flow from singleton paths \\
$m$ & Number of links \\
$n$ & Number of paths (excluding singletons) \\
$s$ & Number of major paths \\
$\ell$ & Number of multi-path OD pairs \\
\midrule
\multicolumn{2}{l}{\textbf{Path Variables}} \\
$x$ & Full path flow vector \\
$y$ & Major-path flow vector \\
$w$ & Minor-path flow vector \\
$z$ & Low-dimensional adjustment coordinate vector \\
\midrule
\multicolumn{2}{l}{\textbf{Incidence Matrices}} \\
$B$ & Path-link incidence matrix ($n \times m$)\\
$B_1$ & Incidence matrix for major paths ($s \times m$)\\
$B_2$ & Incidence matrix for minor paths [$(n-s) \times m$]\\
$A$ & OD-path incidence matrix ($\ell \times n$)\\
$A_1$ & OD-major-path incidence matrix ($\ell \times s$)\\
$A_2$ & OD-minor-path incidence matrix [$\ell \times (n-s)$]\\
$D = B_2'U_r$ & Compressed link matrix ($m \times r$) \\
$M = A_2 U_r$ & Compressed OD matrix ($\ell \times r$) \\
\midrule
\multicolumn{2}{l}{\textbf{Compression Framework}} \\
$\tau$ & Threshold for classifying major and minor paths \\
$B_2 \approx U_r \Sigma_r V_r'$ & Rank-$r$ truncated SVD of minor-path matrix \\
$U_r$ & Matrix of the first $r$ left singular vectors of $B_2$ \\
$r$ & Compression rank \\
$w_0$ & Reference minor-path flow vector \\
$U_rz+w_0$ & Compressed representation of minor-path flows \\
\midrule
\multicolumn{2}{l}{\textbf{Optimization and ALM}} \\
$\hat{f}(y,z)$
& Compressed objective function \\
$\lambda$ & Lagrange multiplier vector (OD equality constraints) \\
$\mu$ & Lagrange multiplier vector associated with $U_rz+w_0\geq0$ \\
$c = (c_1, c_2)$ & Component-wise penalty parameters \\
$\beta$ & Penalty growth factor \\
$\gamma$ & Penalty update threshold \\
\end{longtable}

\section{Gradient Computation and Overhead}
\label{app:gradient-overhead}

The gradient of the AL function with respect to $z$ admits several implementations that trade off memory and computation. We compare four such strategies on the original Chicago {\color{black}Sketch} path pool under nominal demand, using $r=20$ and representative cutoff thresholds from Section~\ref{sec:cutoff-thresholds}.

Under the affine formulation, the terms $B_2'w_0$ in the link flows and $A_2w_0$ in the OD constraints are constant with respect to $z$. Hence, they do not affect the gradient computations below, which involve the same matrices $D=B_2'U_r$ and $M=A_2U_r$.
 
The \emph{direct method} precomputes the matrices $M = A_2 U_r$ (of size $\ell \times r$) and $D = B_2' U_r$ (of size $m \times r$). Although $B_2'$ is sparse, the product $D$ is generally dense, and the per-iteration cost is $O((\ell+m) \cdot r)$.
 
The \emph{chain rule method} evaluates the gradient by composing the link-volume and OD-flow terms without explicitly forming $D$ or $M$. Its computation involves multiplication by $U_r$ together with sparse matrix-vector products involving $B_2$ and $A_2$. Thus, its cost depends on the sparsity of these incidence matrices as well as on $(n-s)r$.
 
The \emph{factored form method} avoids building $D$ explicitly by using the SVD factors directly. Since $ D=B_2'U_r=V_r\Sigma_r$, we have $Dz=V_r\Sigma_r z.$ Thus, the contribution of the low-dimensional adjustment to the link flows can be computed by first scaling $z$ by the singular values and then multiplying by $V_r$. This requires $O(mr)$ operations per evaluation and avoids both the sparse chain through $B_2$ and avoids separately forming \(D\), since \(D=V_r\Sigma_r\) can be applied directly through the stored SVD factors.
 
The \emph{mixed approach}, which we use in our computational studies reported earlier, combines the factored form for the link-volume term with the chain rule for the OD-flow term. This pairing retains the computational efficiency of the factored form on the dominant link-volume computation and while avoiding the additional formation and storage of \(M=A_2U_r\).
 
\begin{table}[htbp]
\centering
\caption{\color{black}Reduction in per-inner-iteration computational time for different gradient implementations on the original Chicago Sketch path pool under nominal demand ($r=20$).}
\label{tab:grad-methods-speedup}
{\small\color{black}
\begin{tabular}{lrrrr}
\toprule
$\tau$ & Direct & Chain Rule & Factored Form & Mixed (Default) \\
\midrule
0.0000 & 0.00\% & 0.00\% & 0.00\% & 0.00\% \\
0.4619 & 27.95\% & 16.98\% & 25.74\% & 25.20\% \\
1.0553 & 35.05\% & 16.18\% & 34.89\% & 30.74\% \\
4.5406 & 35.38\% & 16.89\% & 36.91\% & 33.65\% \\
\bottomrule
\end{tabular}
}
\end{table}

{\color{black}Table~\ref{tab:grad-methods-speedup} reports the reduction in per-inner-iteration computational time of each method relative to the uncompressed baseline. The central observation is that reducing the working dimension from $n-s$ to $r=20$ is considerably more important than preserving sparsity in the original large matrices. At the largest threshold considered, the factored-form and mixed methods reduce the per-inner-iteration computational time by {\color{black}$36.9\%$} and {\color{black}$33.6\%$}, respectively, whereas the chain-rule method achieves a reduction of {\color{black}$16.9\%$}. The direct method performs similarly to the factored form but requires forming and storing the dense matrix $D$.}

\section{Algorithmic Tuning}
\label{app:algorithmic-tuning}

{\color{black}This appendix reports the sensitivity experiments that motivated the default algorithmic settings used throughout Section~\ref{sec:experiments}. All experiments are conducted on the original, non-expanded path pools summarized in Table~\ref{tab:network-characteristics}; the expanded path pools introduced in Section~\ref{sec:path-pool-expansion} are not used here. The major-minor partitions are constructed from the nominal path-flow values, and all reported results use the nominal OD demands rather than perturbed-demand instances.}

{\color{black}In each experiment, we vary a single parameter while holding the others fixed at their default values, and examine the effect on the outer iteration count, total CPU time, and BPR gap. We conduct the sensitivity analysis on Chicago Regional and Philadelphia, the two larger networks on which CPU time is substantial enough for differences between parameter choices to be meaningful. On Chicago Sketch, every configuration we tested solves in under 10 seconds (see Table~\ref{tab:cutoff-sketch} in Appendix~\ref{app:full-results}), and the absolute differences between settings are too small to support reliable conclusions. {\color{black}The effects of the parameter choices differ between networks; the default settings stated in Section~\ref{sec:setup} are retained for consistency across the experiments}.}

\emph{Convergence tolerance:}
We first compare a strict tolerance of $10^{-6}$ with the looser default of $10^{-4}$ on the infinity norm of the constraint violation in the outer iteration. On Chicago Regional at $\tau=1.0329$, the loose tolerance reduces the outer iteration count from {\color{black}13 to 6} and the CPU time from {\color{black}78.31~s to 65.94~s}, a reduction of {\color{black}15.8\%}, while the BPR gap changes by only {\color{black}less than 0.001\%}. On Philadelphia at $\tau=5.3344$, the outer iteration count drops from {\color{black}19 to 10} and the CPU time from {\color{black}170.22~s to 137.02~s}, a reduction of {\color{black}19.5\%}, with a BPR gap change of {\color{black}less than 0.001\%}. Both settings satisfy their prescribed stopping tolerances, indicating that the looser tolerance eliminates additional iterations with only a small effect on solution accuracy.

\emph{Penalty growth factor $\beta$:}
We next compare $\beta=4.0$ and the default $\beta=10.0$. On Chicago Regional at $\tau=1.0329$, the larger growth factor reduces the outer iteration count from {\color{black}8 to 6} {\color{black}but increases} the CPU time from {\color{black}51.93~s to 65.94~s}, corresponding to a {\color{black}27.0\% increase} in computational time, while the BPR gap changes by only {\color{black}0.001\%}. On Philadelphia at $\tau=5.3344$, $\beta=10.0$ reduces the CPU time from {\color{black}162.14~s to 137.02~s}, a reduction of {\color{black}15.5\%}. {\color{black}The effect of the larger growth factor on computational time differs between the two networks.}

\emph{Initial penalty value $c^0$:}
Fixing $\beta=10.0$, we test $c^0_1=c^0_2\in\{10,100,1000\}$ on Chicago Regional at $\tau=1.0329$. The default $c_1^0=c_2^0=1000$ requires {\color{black}6} outer iterations and {\color{black}65.94~s}, whereas $c_1^0=c_2^0=10$ requires {\color{black}13} outer iterations and {\color{black}101.04~s}. The OD constraint violations remain on the order of $10^{-5}$ across all three settings, and the same pattern is observed on Philadelphia. In these experiments, the larger initial penalty provides stronger constraint enforcement from the outset and reduces the computational effort required for convergence.

\emph{Higher inner-iteration budget and penalty-update threshold:}
Finally, we separately examine the effects of increasing the maximum inner iteration count from 200 to 500 and increasing the penalty-update threshold $\gamma$ from 0.25 to 0.50. On Chicago Regional at $\tau=1.0329$, increasing the inner-iteration budget improves the BPR gap from {\color{black}$0.007\%$ to $0.006\%$}, {\color{black}and reduces} CPU time from {\color{black}65.94~s to 51.11~s}, while the link $R^2$ remains essentially unchanged. Raising $\gamma$ to 0.50 instead gives a BPR gap of {\color{black}$0.007\%$} and a CPU time of {\color{black}68.35~s}. {\color{black}Thus, increasing the penalty-update threshold produces only small changes in solution accuracy. Increasing the inner-iteration budget reduces CPU time on Chicago Regional, while its effect differs on Philadelphia. The balance between solution accuracy and computational effort therefore depends on the network.}

\section{Reference Solution Calculations}
\label{app:reference-solution}

The reference solutions are computed using a single-thread gradient-projection procedure on the candidate path pool. Only multi-path OD pairs are optimized, while singleton-path flows are fixed through $v_0$. For the original path pools, the calculations use a uniform cold start within each OD pair. For the X4 and X8 expanded path pools, the calculations start from a converged original-pool assignment: flows on the original paths are retained and the added paths are initialized at zero. Each calculation is terminated when the relative equilibrium gap falls below $10^{-5}$.

\section{Full Cutoff-Threshold Results}
\label{app:full-results}

This appendix reports the complete cutoff-threshold results underlying Section~\ref{sec:cutoff-thresholds}. The experiments use the original path pools under nominal OD demands and the default algorithmic settings with $r=20$.

Across all three networks, substantial reductions in problem dimension are achieved while the computed solutions remain in close numerical agreement with the nominal reference solutions. The BPR gaps remain small, the link-flow $R^2$ values remain close to one, and the travel-time RMSE remains small throughout the tested thresholds. CPU time generally decreases as the compressed dimension becomes smaller, although it is not monotone in $\tau$ because the numbers of outer and inner iterations also vary. These detailed results support the discussion of the effect of the compression threshold in Section~\ref{sec:cutoff-thresholds}.

\begin{table}[htbp]
\centering
\caption{\textcolor{black}{Impact of Cutoff Thresholds (Chicago Sketch, $r=20$).}}
\label{tab:cutoff-sketch}
{\small\color{black}
\resizebox{\textwidth}{!}{%
\begin{tabular}{lrrrrrrrrr}
\toprule
\shortstack{$\tau$ (veh)} & \shortstack{$s$ (major)} & \shortstack{$n-s$ (minor)} & \shortstack{Variable \\ Reduction (\%)} & \shortstack{BPR \\ Gap (\%)} & Link $R^2$ & \shortstack{RMSE (s)} & \shortstack{Outer \\ Iters} & \shortstack{Total Inner \\ Iters} & \shortstack{CPU (s)}  \\
\midrule
0.0000  & 42,774 & 0 & 0.00 & 0.179 & 0.99562 & 4.123 & 7 & 620 & 3.95  \\
0.1385  & 30,149 & 12,625 & 29.47 & 0.039 & 0.99800 & 4.722 & 8 & 854 & 4.45  \\
0.1867  & 29,265 & 13,509 & 31.54 & 0.035 & 0.99812 & 4.409 & 9 & 864 & 4.45  \\
0.2509  & 28,264 & 14,510 & 33.88 & 0.031 & 0.99830 & 3.929 & 6 & 838 & 4.07  \\
0.3405  & 27,050 & 15,724 & 36.71 & 0.030 & 0.99832 & 3.889 & 6 & 622 & 2.97  \\
0.4619  & 25,634 & 17,140 & 40.02 & 0.029 & 0.99834 & 3.833 & 6 & 720 & 3.16  \\
0.6766  & 24,185 & 18,589 & 43.41 & 0.029 & 0.99835 & 3.645 & 7 & 842 & 3.56  \\
1.0553  & 22,714 & 20,060 & 46.85 & 0.028 & 0.99843 & 3.544 & 8 & 721 & 2.99  \\
1.8918  & 21,068 & 21,706 & 50.70 & 0.025 & 0.99856 & 3.383 & 9 & 564 & 2.35  \\
4.5406  & 19,347 & 23,427 & 54.72 & 0.027 & 0.99854 & 3.952 & 7 & 599 & 2.28  \\
\bottomrule
\end{tabular}%
}}
\end{table}

\begin{table}[htbp]
\centering
\caption{\textcolor{black}{Impact of Cutoff Thresholds (Chicago Regional, $r=20$).}}
\label{tab:cutoff-regional}
{\small\color{black}
\resizebox{\textwidth}{!}{%
\begin{tabular}{lrrrrrrrrr}
\toprule
\shortstack{$\tau$ (veh)} & \shortstack{$s$ (major)} & \shortstack{$n-s$ (minor)} & \shortstack{Variable \\ Reduction (\%)} & \shortstack{BPR \\ Gap (\%)} & Link $R^2$ & \shortstack{RMSE (s)} & \shortstack{Outer \\ Iters} & \shortstack{Total Inner \\ Iters} & \shortstack{CPU (s)}  \\
\midrule
0.0000  & 879,625 & 0 & 0.00 & 0.102 & 0.99939 & 1.653 & 6 & 817 & 111.28  \\
0.1264  & 546,723 & 332,902 & 37.84 & 0.008 & 0.99990 & 1.118 & 7 & 714 & 81.24  \\
0.1535  & 531,721 & 347,904 & 39.55 & 0.008 & 0.99990 & 1.148 & 8 & 681 & 77.27  \\
0.1857  & 514,603 & 365,022 & 41.50 & 0.008 & 0.99990 & 1.058 & 7 & 741 & 81.38  \\
0.2265  & 495,323 & 384,302 & 43.69 & 0.009 & 0.99989 & 1.139 & 7 & 454 & 51.43  \\
0.2798  & 474,886 & 404,739 & 46.01 & 0.008 & 0.99991 & 1.099 & 7 & 639 & 67.32  \\
0.3534  & 453,041 & 426,584 & 48.49 & 0.007 & 0.99991 & 1.048 & 6 & 619 & 62.90  \\
0.4602  & 430,281 & 449,344 & 51.08 & 0.007 & 0.99992 & 1.098 & 6 & 616 & 61.77  \\
0.6367  & 407,221 & 472,404 & 53.70 & 0.007 & 0.99992 & 1.057 & 6 & 556 & 54.94  \\
1.0329  & 382,377 & 497,248 & 56.53 & 0.007 & 0.99993 & 0.986 & 6 & 691 & 65.94  \\
\bottomrule
\end{tabular}%
}}
\end{table}

\begin{table}[htbp]
\centering
\caption{\textcolor{black}{Impact of Cutoff Thresholds (Philadelphia, $r=20$).}}
\label{tab:cutoff-philly}
{\small\color{black}
\resizebox{\textwidth}{!}{%
\begin{tabular}{lrrrrrrrrr}
\toprule
\shortstack{$\tau$ (veh)} & \shortstack{$s$ (major)} & \shortstack{$n-s$ (minor)} & \shortstack{Variable \\ Reduction (\%)} & \shortstack{BPR \\ Gap (\%)} & Link $R^2$ & \shortstack{RMSE (s)} & \shortstack{Outer \\ Iters} & \shortstack{Total Inner \\ Iters} & \shortstack{CPU (s)}  \\
\midrule
0.0000  & 1,846,068 & 0 & 0.00 & 0.414 & 0.99559 & 4.717 & 13 & 1,034 & 315.55  \\
0.1610  & 1,100,488 & 745,580 & 40.39 & 0.065 & 0.99907 & 3.328 & 8 & 1,030 & 231.89  \\
0.2239  & 1,060,985 & 785,083 & 42.53 & 0.064 & 0.99909 & 3.310 & 6 & 837 & 185.22  \\
0.3117  & 1,014,000 & 832,068 & 45.07 & 0.071 & 0.99903 & 3.523 & 7 & 953 & 204.63  \\
0.4473  & 956,142 & 889,926 & 48.21 & 0.067 & 0.99907 & 3.429 & 9 & 829 & 177.86  \\
0.6663  & 899,580 & 946,488 & 51.27 & 0.061 & 0.99915 & 3.252 & 8 & 851 & 177.08  \\
0.9866  & 840,996 & 1,005,072 & 54.44 & 0.061 & 0.99916 & 3.252 & 8 & 842 & 168.69  \\
1.4643  & 787,486 & 1,058,582 & 57.34 & 0.064 & 0.99916 & 3.351 & 8 & 943 & 182.89  \\
2.4435  & 731,407 & 1,114,661 & 60.38 & 0.054 & 0.99929 & 3.052 & 8 & 690 & 132.42  \\
5.3344  & 671,735 & 1,174,333 & 63.61 & 0.054 & 0.99935 & 3.101 & 10 & 718 & 137.02  \\
\bottomrule
\end{tabular}%
}}
\end{table}
 
\section{Path and Flow Distribution Tables}
\label{app:pathflow}
 
This appendix reports the nominal flows carried by the major and minor paths at selected cutoff thresholds, together with their shares of total demand, complementing the iteration and CPU results of Appendix~\ref{app:full-results}. The three networks display markedly different concentration profiles.

On Chicago Sketch (Table~\ref{tab:pathflow-sketch}), major paths carry roughly 23\% of total demand across all thresholds and minor-path flow stays below {\color{black}0.9\%} even at the most aggressive cutoff, indicating that the vast majority of demand is absorbed by singleton OD pairs that lie outside the compressed problem. On Chicago Regional (Table~\ref{tab:pathflow-regional}), major paths carry {\color{black}54--59\%} of demand, with minor-path flow rising to {\color{black}5.0\%} at $\tau = 1.033$. On Philadelphia (Table~\ref{tab:pathflow-philly}), major paths dominate at {\color{black}67--71\%} of demand and minor-path flow reaches {\color{black}3.6\%} at $\tau = 5.3344$. {At the largest cutoffs, the variable reductions are 56.53\% on Chicago Regional and 63.61\% on Philadelphia, whereas the corresponding minor-flow shares are 5.0\% and 3.6\%. Thus, variable reduction and the minor-path flow share describe different properties of the major--minor partition.} In all three networks, the fraction of total flow assigned to minor paths increases progressively with $\tau$, illustrating how the cutoff controls both the number of paths assigned to the minor-path block and their share of the nominal flow.
 
\begin{table}[htbp]
\centering
\caption{Path and Nominal-Flow Distribution (Chicago Sketch)}
\label{tab:pathflow-sketch}
{\small\color{black}
\begin{tabular}{rrrrr}
\toprule
\shortstack{$\tau$ (veh)} & \shortstack{$s$ (major)} & \shortstack{$n-s$ (minor)} & Major Flow (\%) & Minor Flow (\%) \\
\midrule
0.000 & 42{,}774 & 0 & 267{,}805.15 (23.5) & 0.00 (0.0) \\
0.139 & 30{,}149 & 12{,}625 & 267{,}683.26 (23.5) & 121.89 (0.0) \\
0.341 & 27{,}050 & 15{,}724 & 266{,}963.93 (23.5) & 841.22 (0.1) \\
1.055 & 22{,}714 & 20{,}060 & 264{,}342.67 (23.2) & 3{,}462.49 (0.3) \\
4.541 & 19{,}347 & 23{,}427 & 257{,}025.34 (22.6) & 10{,}779.81 (0.9) \\
\bottomrule
\end{tabular}}
\end{table}
 
\begin{table}[htbp]
\centering
\caption{Path and Nominal-Flow Distribution (Chicago Regional)}
\label{tab:pathflow-regional}
{\small\color{black}
\begin{tabular}{rrrrr}
\toprule
\shortstack{$\tau$ (veh)} & \shortstack{$s$ (major)} & \shortstack{$n-s$ (minor)} & Major Flow (\%) & Minor Flow (\%) \\
\midrule
0.000 & 879{,}625 & 0 & 778{,}822.96 (59.2) & 0.00 (0.0) \\
0.126 & 546{,}723 & 332{,}902 & 775{,}629.74 (58.9) & 3{,}193.22 (0.2) \\
0.227 & 495{,}323 & 384{,}302 & 766{,}665.31 (58.3) & 12{,}157.65 (0.9) \\
0.460 & 430{,}281 & 449{,}344 & 745{,}455.32 (56.6) & 33{,}367.64 (2.5) \\
1.033 & 382{,}377 & 497{,}248 & 713{,}025.41 (54.2) & 65{,}797.55 (5.0) \\
\bottomrule
\end{tabular}}
\end{table}
 
\begin{table}[htbp]
\centering
\caption{Path and Nominal-Flow Distribution (Philadelphia)}
\label{tab:pathflow-philly}
{\small\color{black}
\begin{tabular}{rrrrr}
\toprule
\shortstack{$\tau$ (veh)} & \shortstack{$s$ (major)} & \shortstack{$n-s$ (minor)} & Major Flow (\%) & Minor Flow (\%) \\
\midrule
0.0000 & 1{,}846{,}068 & 0 & 10{,}157{,}631.02 (70.9) & 0.00 (0.0) \\
0.1610 & 1{,}100{,}488 & 745{,}580 & 10{,}151{,}617.68 (70.9) & 6{,}013.34 (0.0) \\
0.4473 & 956{,}142 & 889{,}926 & 10{,}109{,}710.00 (70.6) & 47{,}921.02 (0.3) \\
0.9866 & 840{,}996 & 1{,}005{,}072 & 10{,}030{,}939.11 (70.0) & 126{,}691.91 (0.9) \\
2.4435 & 731{,}407 & 1{,}114{,}661 & 9{,}859{,}793.75 (68.8) & 297{,}837.27 (2.1) \\
5.3344 & 671{,}735 & 1{,}174{,}333 & 9{,}645{,}316.17 (67.3) & 512{,}314.86 (3.6) \\
\bottomrule
\end{tabular}}
\end{table}

\section{{Full Demand Perturbation Results}}
\label{app:demand-perturbation-full}

{This appendix reports the complete seed-level results underlying Tables~\ref{tab:demand-perturbations-original} and \ref{tab:demand-perturbations-expansion}. The notation ``Uncomp./Comp.'' gives the uncompressed and compressed values, respectively. CPU reduction is the percentage reduction in computational time of the compressed formulation relative to the corresponding uncompressed solve. For each perturbation experiment, the major--minor partition and rank-$20$ SVD basis are constructed from the nominal-demand problem and kept fixed after the OD demands are perturbed, while the reference minor-path flows $w_0$ are scaled by the same factors as the corresponding OD demands.}

\begin{table}[htbp]
\centering
\caption{\textcolor{black}{Full demand perturbation results on the original path pools ($r=20$).}}
\label{tab:demand-perturbations-original-full}
{\scriptsize\color{black}
\resizebox{\textwidth}{!}{%
\begin{tabular}{lrrrrrrrrr}
\toprule
Network & Seed & Pert. & \shortstack{Variable\\Reduction (\%)} & \shortstack{BPR Gap (\%)\\Uncomp. / Comp.} & \shortstack{Link $R^2$\\Uncomp. / Comp.} & \shortstack{RMSE (s)\\Uncomp. / Comp.} & \shortstack{Outer / Inner\\Uncomp. / Comp.} & \shortstack{CPU (s)\\Uncomp. / Comp.} & \shortstack{CPU\\Reduction (\%)}  \\
\midrule
Chicago Sketch & 10 & $\pm5\%$  & 54.72 & $0.001 / 0.001$ & $0.99998 / 0.99998$ & $0.332 / 0.383$ & $4/409\;/\;4/195$ & $2.37 / 0.82$ & 65.37  \\
 & 28 & $\pm5\%$  & 54.72 & $0.001 / 0.001$ & $0.99998 / 0.99998$ & $0.275 / 0.293$ & $4/409\;/\;6/158$ & $2.38 / 0.73$ & 69.34  \\
 & 42 & $\pm5\%$  & 54.72 & $0.001 / 0.001$ & $0.99998 / 0.99998$ & $0.338 / 0.309$ & $4/409\;/\;6/162$ & $2.42 / 0.76$ & 68.53  \\
 & 10 & $\pm10\%$  & 54.72 & $0.001 / 0.001$ & $0.99998 / 0.99997$ & $0.436 / 0.673$ & $12/429\;/\;6/182$ & $3.05 / 0.83$ & 72.67  \\
 & 28 & $\pm10\%$  & 54.72 & $0.001 / 0.001$ & $0.99998 / 0.99997$ & $0.358 / 0.450$ & $4/409\;/\;6/182$ & $2.36 / 0.84$ & 64.42  \\
 & 42 & $\pm10\%$  & 54.72 & $0.001 / 0.001$ & $0.99998 / 0.99997$ & $0.385 / 0.488$ & $7/416\;/\;5/159$ & $2.64 / 0.71$ & 73.10  \\
Chicago Regional & 10 & $\pm5\%$  & 56.53 & $0.002 / 0.002$ & $0.99998 / 0.99998$ & $0.308 / 0.415$ & $6/617\;/\;4/354$ & $82.76 / 34.54$ & 58.27  \\
 & 28 & $\pm5\%$  & 56.53 & $0.002 / 0.002$ & $0.99998 / 0.99998$ & $0.310 / 0.462$ & $6/617\;/\;5/350$ & $82.70 / 34.66$ & 58.09  \\
 & 42 & $\pm5\%$  & 56.53 & $0.002 / 0.002$ & $0.99998 / 0.99998$ & $0.301 / 0.443$ & $8/621\;/\;5/352$ & $86.17 / 35.29$ & 59.05  \\
 & 10 & $\pm10\%$  & 56.53 & $0.002 / 0.002$ & $0.99998 / 0.99999$ & $0.317 / 0.400$ & $9/624\;/\;5/202$ & $88.34 / 21.92$ & 75.19  \\
 & 28 & $\pm10\%$  & 56.53 & $0.002 / 0.002$ & $0.99998 / 0.99998$ & $0.337 / 0.406$ & $6/618\;/\;6/344$ & $82.71 / 34.75$ & 57.99  \\
 & 42 & $\pm10\%$  & 56.53 & $0.002 / 0.002$ & $0.99998 / 0.99999$ & $0.306 / 0.490$ & $7/619\;/\;9/264$ & $84.71 / 30.55$ & 63.94  \\
Philadelphia & 10 & $\pm5\%$  & 63.61 & $0.002 / 0.002$ & $0.99999 / 0.99999$ & $0.278 / 0.317$ & $6/616\;/\;6/174$ & $174.40 / 37.90$ & 78.27  \\
 & 28 & $\pm5\%$  & 63.61 & $0.002 / 0.002$ & $0.99999 / 0.99999$ & $0.292 / 0.341$ & $10/625\;/\;8/186$ & $190.17 / 43.09$ & 77.34  \\
 & 42 & $\pm5\%$  & 63.61 & $0.002 / 0.002$ & $0.99999 / 0.99999$ & $0.276 / 0.311$ & $7/615\;/\;8/181$ & $179.39 / 43.17$ & 75.94  \\
 & 10 & $\pm10\%$  & 63.61 & $0.002 / 0.002$ & $0.99998 / 0.99998$ & $0.495 / 0.488$ & $7/425\;/\;5/271$ & $130.21 / 53.78$ & 58.70  \\
 & 28 & $\pm10\%$  & 63.61 & $0.002 / 0.002$ & $0.99998 / 0.99998$ & $0.444 / 0.460$ & $6/614\;/\;6/351$ & $173.75 / 67.00$ & 61.44  \\
 & 42 & $\pm10\%$  & 63.61 & $0.002 / 0.003$ & $0.99998 / 0.99998$ & $0.421 / 0.503$ & $7/620\;/\;5/256$ & $178.25 / 50.42$ & 71.72  \\
\bottomrule
\end{tabular}%
}}
\end{table}

\begin{table}[htbp]
\centering
\caption{\textcolor{black}{Full demand perturbation results on the expanded X4 and X8 path pools for Chicago Sketch ($r=20$).}}
\label{tab:demand-perturbations-expansion-sketch-full}
{\scriptsize\color{black}
\resizebox{\textwidth}{!}{%
\begin{tabular}{lrrrrrrrrr}
\toprule
Network & Seed & Pert. & \shortstack{Variable\\Reduction (\%)} & \shortstack{BPR Gap (\%)\\Uncomp. / Comp.} & \shortstack{Link $R^2$\\Uncomp. / Comp.} & \shortstack{RMSE (s)\\Uncomp. / Comp.} & \shortstack{Outer / Inner\\Uncomp. / Comp.} & \shortstack{CPU (s)\\Uncomp. / Comp.} & \shortstack{CPU\\Reduction (\%)}  \\
\midrule
Chicago Sketch X4 & 10 & $\pm5\%$  & 63.81 & $0.002 / 0.002$ & $0.99997 / 0.99997$ & $0.355 / 0.469$ & $4/410\;/\;4/122$ & $3.36 / 0.80$ & 76.04  \\
 & 28 & $\pm5\%$  & 63.81 & $0.002 / 0.002$ & $0.99997 / 0.99997$ & $0.351 / 0.324$ & $4/410\;/\;4/137$ & $3.32 / 0.87$ & 73.72  \\
 & 42 & $\pm5\%$  & 63.81 & $0.002 / 0.002$ & $0.99997 / 0.99997$ & $0.385 / 0.393$ & $4/410\;/\;6/163$ & $3.35 / 1.11$ & 66.70  \\
 & 10 & $\pm10\%$  & 63.81 & $0.002 / 0.002$ & $0.99997 / 0.99997$ & $0.467 / 0.522$ & $4/410\;/\;7/168$ & $3.35 / 1.18$ & 64.79  \\
 & 28 & $\pm10\%$  & 63.81 & $0.002 / 0.002$ & $0.99997 / 0.99997$ & $0.407 / 0.689$ & $4/410\;/\;5/174$ & $3.39 / 1.18$ & 65.16  \\
 & 42 & $\pm10\%$  & 63.81 & $0.002 / 0.002$ & $0.99997 / 0.99997$ & $0.458 / 0.833$ & $4/410\;/\;7/126$ & $3.38 / 1.01$ & 70.01  \\
Chicago Sketch X8 & 10 & $\pm5\%$  & 81.67 & $0.002 / 0.002$ & $0.99997 / 0.99997$ & $0.382 / 0.867$ & $4/415\;/\;5/131$ & $5.71 / 1.63$ & 71.55  \\
 & 28 & $\pm5\%$  & 81.67 & $0.002 / 0.002$ & $0.99997 / 0.99997$ & $0.362 / 0.745$ & $4/415\;/\;4/97$ & $5.72 / 1.22$ & 78.70  \\
 & 42 & $\pm5\%$  & 81.67 & $0.002 / 0.002$ & $0.99997 / 0.99997$ & $0.396 / 0.418$ & $4/415\;/\;7/92$ & $5.72 / 1.32$ & 77.01  \\
 & 10 & $\pm10\%$  & 81.67 & $0.002 / 0.002$ & $0.99996 / 0.99997$ & $0.496 / 0.889$ & $4/415\;/\;4/234$ & $5.77 / 2.59$ & 55.02  \\
 & 28 & $\pm10\%$  & 81.67 & $0.002 / 0.002$ & $0.99997 / 0.99996$ & $0.434 / 0.580$ & $4/415\;/\;5/80$ & $5.73 / 1.05$ & 81.62  \\
 & 42 & $\pm10\%$  & 81.67 & $0.002 / 0.002$ & $0.99996 / 0.99996$ & $0.480 / 0.449$ & $4/415\;/\;5/77$ & $5.72 / 1.03$ & 81.94  \\
\bottomrule
\end{tabular}%
}}
\end{table}

\begin{table}[htbp]
\centering
\caption{\textcolor{black}{Full X4/X8 demand perturbation results on Chicago Regional ($r=20$).}}
\label{tab:demand-perturbations-expansion-regional-full}
{\scriptsize\color{black}
\resizebox{\textwidth}{!}{%
\begin{tabular}{lrrrrrrrrr}
\toprule
Network & Seed & Pert. & \shortstack{Variable\\Reduction (\%)} & \shortstack{BPR Gap (\%)\\Uncomp. / Comp.} & \shortstack{Link $R^2$\\Uncomp. / Comp.} & \shortstack{RMSE (s)\\Uncomp. / Comp.} & \shortstack{Outer / Inner\\Uncomp. / Comp.} & \shortstack{CPU (s)\\Uncomp. / Comp.} & \shortstack{CPU\\Reduction (\%)}  \\
\midrule
Chicago Regional X4 & 10 & $\pm5\%$  & 65.73 & $0.003 / 0.003$ & $0.99998 / 0.99998$ & $0.341 / 0.305$ & $6/812\;/\;5/528$ & $155.57 / 72.64$ & 53.31  \\
 & 28 & $\pm5\%$  & 65.73 & $0.003 / 0.002$ & $0.99998 / 0.99998$ & $0.340 / 0.383$ & $6/812\;/\;6/616$ & $155.62 / 85.40$ & 45.12  \\
 & 42 & $\pm5\%$  & 65.73 & $0.003 / 0.003$ & $0.99998 / 0.99998$ & $0.347 / 0.459$ & $6/813\;/\;5/400$ & $155.21 / 56.49$ & 63.60  \\
 & 10 & $\pm10\%$  & 65.73 & $0.003 / 0.002$ & $0.99998 / 0.99998$ & $0.355 / 0.319$ & $6/812\;/\;9/624$ & $154.92 / 88.80$ & 42.68  \\
 & 28 & $\pm10\%$  & 65.73 & $0.003 / 0.003$ & $0.99998 / 0.99998$ & $0.353 / 0.365$ & $6/812\;/\;6/484$ & $153.92 / 67.05$ & 56.44  \\
 & 42 & $\pm10\%$  & 65.73 & $0.003 / 0.003$ & $0.99998 / 0.99998$ & $0.352 / 0.373$ & $6/812\;/\;4/428$ & $153.87 / 58.57$ & 61.93  \\
Chicago Regional X8 & 10 & $\pm5\%$  & 82.63 & $0.004 / 0.003$ & $0.99998 / 0.99998$ & $0.354 / 0.371$ & $6/813\;/\;6/620$ & $279.02 / 137.71$ & 50.64  \\
 & 28 & $\pm5\%$  & 82.63 & $0.004 / 0.003$ & $0.99998 / 0.99998$ & $0.354 / 0.356$ & $6/813\;/\;6/618$ & $278.88 / 135.94$ & 51.26  \\
 & 42 & $\pm5\%$  & 82.63 & $0.004 / 0.003$ & $0.99998 / 0.99998$ & $0.348 / 0.356$ & $6/814\;/\;5/623$ & $283.18 / 136.09$ & 51.94  \\
 & 10 & $\pm10\%$  & 82.63 & $0.004 / 0.003$ & $0.99997 / 0.99998$ & $0.370 / 0.419$ & $6/814\;/\;6/618$ & $279.38 / 135.82$ & 51.39  \\
 & 28 & $\pm10\%$  & 82.63 & $0.004 / 0.003$ & $0.99997 / 0.99998$ & $0.375 / 0.402$ & $6/813\;/\;5/580$ & $278.17 / 126.89$ & 54.38  \\
 & 42 & $\pm10\%$  & 82.63 & $0.004 / 0.003$ & $0.99997 / 0.99998$ & $0.357 / 0.348$ & $6/813\;/\;7/620$ & $282.63 / 138.38$ & 51.04  \\
\bottomrule
\end{tabular}%
}}
\end{table}

\begin{table}[htbp]
\centering
\caption{\textcolor{black}{Full X4/X8 demand perturbation results on Philadelphia ($r=20$).}}
\label{tab:demand-perturbations-expansion-philadelphia-full}
{\scriptsize\color{black}
\resizebox{\textwidth}{!}{%
\begin{tabular}{lrrrrrrrrr}
\toprule
Network & Seed & Pert. & \shortstack{Variable\\Reduction (\%)} & \shortstack{BPR Gap (\%)\\Uncomp. / Comp.} & \shortstack{Link $R^2$\\Uncomp. / Comp.} & \shortstack{RMSE (s)\\Uncomp. / Comp.} & \shortstack{Outer / Inner\\Uncomp. / Comp.} & \shortstack{CPU (s)\\Uncomp. / Comp.} & \shortstack{CPU\\Reduction (\%)}  \\
\midrule
Philadelphia X4 & 10 & $\pm5\%$  & 65.38 & $0.004 / 0.004$ & $0.99994 / 0.99994$ & $0.383 / 0.439$ & $7/623\;/\;6/461$ & $255.40 / 117.65$ & 53.93  \\
 & 28 & $\pm5\%$  & 65.38 & $0.004 / 0.004$ & $0.99994 / 0.99994$ & $0.386 / 0.456$ & $6/618\;/\;4/429$ & $248.75 / 107.40$ & 56.82  \\
 & 42 & $\pm5\%$  & 65.38 & $0.004 / 0.005$ & $0.99994 / 0.99994$ & $0.379 / 0.543$ & $6/618\;/\;9/384$ & $249.14 / 104.92$ & 57.89  \\
 & 10 & $\pm10\%$  & 65.38 & $0.005 / 0.005$ & $0.99993 / 0.99993$ & $0.498 / 0.627$ & $7/623\;/\;4/360$ & $255.16 / 92.22$ & 63.86  \\
 & 28 & $\pm10\%$  & 65.38 & $0.005 / 0.005$ & $0.99992 / 0.99992$ & $0.516 / 0.578$ & $6/619\;/\;6/444$ & $251.74 / 115.01$ & 54.31  \\
 & 42 & $\pm10\%$  & 65.38 & $0.005 / 0.005$ & $0.99993 / 0.99993$ & $0.491 / 0.644$ & $7/621\;/\;6/379$ & $253.72 / 98.45$ & 61.20  \\
Philadelphia X8 & 10 & $\pm5\%$  & 81.40 & $0.005 / 0.005$ & $0.99992 / 0.99992$ & $0.401 / 0.615$ & $6/623\;/\;8/428$ & $480.02 / 178.77$ & 62.76  \\
 & 28 & $\pm5\%$  & 81.40 & $0.005 / 0.005$ & $0.99992 / 0.99992$ & $0.392 / 0.592$ & $6/623\;/\;8/432$ & $472.29 / 175.90$ & 62.76  \\
 & 42 & $\pm5\%$  & 81.40 & $0.005 / 0.005$ & $0.99992 / 0.99992$ & $0.390 / 0.599$ & $6/622\;/\;6/426$ & $467.65 / 168.82$ & 63.90  \\
 & 10 & $\pm10\%$  & 81.40 & $0.005 / 0.006$ & $0.99991 / 0.99991$ & $0.553 / 0.755$ & $6/623\;/\;6/437$ & $475.74 / 174.74$ & 63.27  \\
 & 28 & $\pm10\%$  & 81.40 & $0.005 / 0.006$ & $0.99991 / 0.99991$ & $0.532 / 0.692$ & $6/623\;/\;8/429$ & $501.05 / 181.56$ & 63.76  \\
 & 42 & $\pm10\%$  & 81.40 & $0.005 / 0.006$ & $0.99991 / 0.99991$ & $0.522 / 0.701$ & $6/623\;/\;7/441$ & $484.65 / 179.38$ & 62.99  \\
\bottomrule
\end{tabular}%
}}
\end{table}

\end{document}